\documentclass[11pt]{article}
\usepackage{amsfonts}
\usepackage[utf8]{inputenc}

\usepackage{graphics}
\usepackage{color}
\usepackage{indentfirst}
\usepackage{cite}
\usepackage{latexsym}
\usepackage{amsmath}
\usepackage{amssymb}
\usepackage[dvips]{epsfig}
\usepackage{amscd}

\newtheorem{theorem}{Theorem}[section]
\newtheorem{remark}{Remark}[section]
\newcommand{\thatsall}{\hfill$\Box$}
\newtheorem{definition}{Definition}[section]
\newtheorem{lemma}[theorem]{Lemma}

\newtheorem{corollary}[theorem]{Corollary}
\newtheorem{proposition}[theorem]{Proposition}

\newcommand{\n}{\rho}

\newcommand{\ti}{\tilde}

\newcommand{\lm}{\lambda}

\renewcommand{\div}{ {\rm div }  }

\newcommand{\pa}{\partial}
\renewcommand{\r}{\mathbb{R}}
\newcommand{\bi}{\bibitem}

\newcommand{\ia}{\int_0^T}

\newcommand{\bt}{\begin{theorem}}
\newcommand{\bl}{\begin{lemma}}
\newcommand{\el}{\end{lemma}}
\newcommand{\et}{\end{theorem}}
\newcommand{\ga}{\gamma}

\newcommand{\curl}{{\rm curl} }

\newcommand{\al}{\alpha}
\newcommand{\de}{\delta}
\newcommand{\ve}{\varepsilon}
\newcommand{\la}{\label}

\newcommand{\ol}{\overline}

\newcommand{\bn}{\begin{eqnarray}}
\newcommand{\en}{\end{eqnarray}}
\newcommand{\bnn}{\begin{eqnarray*}}
\newcommand{\enn}{\end{eqnarray*}}

\newcommand{\bnnn}{\begin{eqnarray*}}
\newcommand{\ennn}{\end{eqnarray*}}

\newcommand{\ba}{\begin{aligned}}
\newcommand{\ea}{\end{aligned}}
\newcommand{\be}{\begin{equation}}
\newcommand{\ee}{\end{equation}}
\def\O{{\Omega }}

\def\norm[#1]#2{\|#2\|_{#1}}

\def\o{\omega}

\newcommand{\si}{\sigma}

\def\la{\label}

\def\na{\nabla}

\makeatletter      
\@addtoreset{equation}{section}
\makeatother

\makeatletter      
\@addtoreset{equation}{section}
\makeatother       

\title{ Existence and Exponential Growth of  Global Classical Solutions to the  Compressible Navier-Stokes Equations with Slip Boundary Conditions in 3D Bounded Domains}

 \author{Guocai C{\small AI}$^{a}$,
   Jing L{\small I}$^{ b,c,d} $ \thanks{email:  gotry@xmu.edu.cn (G. C. Cai), ajingli@gmail.com (J. Li)}  \\
{\normalsize a.  School of Mathematical Sciences, }\\ {\normalsize  Xiamen University, Xiamen 361005, P. R. China;}\\
{\normalsize b. Department of Mathematics, }\\ {\normalsize  Nanchang University, Nanchang 330031, P. R. China;} \\ {\normalsize c. Institute of Applied Mathematics, AMSS,} \\ {\normalsize \&   Hua Loo-Keng Key Laboratory of Mathematics,}\\
{\normalsize  Chinese Academy of Sciences,    Beijing 100190,
P. R. China;}
 \\ {\normalsize d.  School of Mathematical Sciences,}\\
{\normalsize  University of Chinese Academy of Sciences, Beijing 100049, P. R. China}}

\date{ }

\begin{document}
\maketitle

 \begin{abstract}
 We investigate  the barotropic compressible Navier-Stokes equations  with  slip boundary conditions in a three-dimensional (3D) bounded domain, whose smooth boundary has a finite number of two-dimensional connected components. For any   adiabatic exponent bigger than one, after  discovering some new estimates  on     boundary integrals related to the  slip boundary condition, we prove that both the weak and classical solutions to the initial-boundary-value problem of this system   exist  globally in time provided the initial energy is suitably small. Moreover, the density  has large oscillations and contains  vacuum states. Finally,  it is also shown that for the classical solutions, the oscillation of the density will grow unboundedly in the long run with an exponential  rate provided  vacuum appears  (even at a point) initially. This is the first result concerning the global existence   of   classical solutions to the compressible Navier-Stokes equations   with density containing vacuum  states initially for  general 3D bounded smooth domains.
 \end{abstract}


\section{Introduction}

The viscous barotropic compressible Navier-Stokes equations for isentropic flows express the principles of conservation of mass and momentum in the absence of exterior forces:
   \be \la{a1}  \begin{cases}\n_t+{\rm div} (\n u)=0,\\
 (\n u)_t+{\rm div}(\n u\otimes u)-\mu\Delta u-(\mu+\lambda)\na {\rm
div} u +\na P(\n) =0, \end{cases}\ee
 where $(x,t)\in\Omega\times (0,T]$, $\Omega$ is a domain in $\r^{N}$, $t\ge 0$ is time,    $x$  is the spatial coordinate, and $\rho\geq0,\,\, u=(u^1,\cdots,u^N)$
are the unknown fluid density, velocity respectively.
 The constants $\mu$ and $\lambda$ are the shear and bulk viscosity coefficients respectively  satisfying the following physical restrictions: \be\la{h3} \mu>0,\quad 2\mu +  {N} \lambda\ge 0.
\ee
 We consider the barotropic case, that is, the pressure $P(\n)$ satisfies \be\la{h4} P(\rho)=a\rho^{\gamma}\ee with constants  $a>0$ and  $\gamma>1 $ which is the adiabatic exponent.

In this paper, we  assume that $\Omega$ is a bounded domain in $\r^3$, its boundary $\partial\Omega$ is of class $C^{\infty}$ and only has a finite number of 2-dimensional connected components. In addition, the system is studied subject to the given initial data
\be \la{h2} \n(x,0)=\n_0(x), \quad \n u(x,0)= \n_0u_0(x),\quad x\in \Omega,\ee
and slip boundary condition
\be \la{ch1} u \cdot n = 0, \,\,\curl u\times n=-A\, u \,\,\,\text{on}\,\,\, \partial\Omega,\ee
where $A=A(x)$ is a $3\times 3$ symmetric matrix defined on $\partial\Omega$.

The mathematical study of compressible Navier-Stokes equations dates back
to the late 1950s.  For   density away from vacuum, Serrin \cite{se1} and Nash \cite{Na}  first  considered
the mathematical questions of compressible viscous fluid dynamics. An intensive
treatment of compressible Navier-Stokes equations started with pioneering works
by Itaya \cite{nit1}, Matsumura-Nishida \cite{M1}, Kazhikhov-Solonnikov \cite{Kaz2}, and Hoff \cite{Hof2} on the local theory for nonstationary problems, and by Beir\~{a}o da Veiga \cite{hbdv1,hbdv2}, Padula \cite{MPa1}, and Novotn\'{y}-Padula \cite{AM1,AM2} on the theory of stationary problems for small data. For the case that  density contains vacuum,  
 Lions \cite{L1}   proved   the global existence of so called finite-energy weak solutions when the adiabatic exponent $\gamma$ is suitably large, for example, $\gamma \geq  {9}/{5}$ for 3D case.   These results were   further improved   by Feireisl-Novotn\'{y}-Petzeltov\'{a}  \cite{fe1} to $\ga>3/2 $  for  three-dimensional case.  Moreover,  Hoff \cite{H3,Hof2,Ho3,HS1,HT1} considered a new type of global weak solutions for any  $\ga>1 $  with small energy that have extra regularity information compared with Lions-Feireisl's large weak ones. However, the regularity and
uniqueness of those weak solutions  \cite{L1, fe1,Ho3,HS1,HT1}  are completely open.  
Recently,  Huang-Li-Xin \cite{hlx1} and Li-Xin \cite{jx01} established the global well-posedness
of classical solutions to the Cauchy problem for the 3D and 2D barotropic compressible Navier-Stokes equations in whole space with smooth initial data that are of small energy but
possibly large oscillations, in particular, the initial density is allowed to vanish. However, when   the domains are bounded, the global  existence of  strong and/or classical solutions with vacuum to the compressible Navier-Stokes equations  \eqref{a1} remains open.

For bounded domains, the usual  Navier-type slip condition  can be stated as follows:
\be \la{Navi}
u \cdot n = 0, \,\,(2D(u)\,n+ \vartheta u)_{\rm tan }=0 \,\,\,\text{on}\,\,\, \partial\Omega,
\ee
where $D(u) = (\nabla u+(\nabla u)^{\rm tr})/2$ is the shear stress, $\vartheta$ is a scalar friction function  which measures the tendency of the fluid to slip on the boundary, and the symbol $v_{\rm tan }$ represents the projection of tangent plane of the vector $v$ on $\partial\Omega$.  Indeed,  introduced originally by Navier \cite{Nclm1}   and later independently by Maxwell \cite{MJC}, the   Navier-type slip condition \eqref{Navi}, which shows that there is a stagnant layer of fluid close to the wall allowing a fluid to slip with the slip velocity being proportional to the shear stress,    is frequently used in numerical studies and analysis for various fluid mechanical problems ( see  for instance \cite{Cpfc1, Itt2, se2} and their references therein).
Especially to deserve to be mentioned, the restriction that $\vartheta$ is non negative
is usual, in order to ensure the conservation of energy. But mathematically, we can take into
account the negative values of $\vartheta$ as well.
For compressible Navier-Stokes equations \eqref{a1} with Navier-type slip boundary condition \eqref{Navi}, Novoton\'y-Stra\v{s}kraba \cite{ANIS} obtained the the global  existence of weak solutions for   $\ga>3/2$ and $\vartheta=0$ in a non axis-symmetric domain when $\mu>0$ and $2\mu +  3 \lambda> 0$.  Hoff \cite{Ho3} studied the global existence of weak solutions on the half space in $\r^3$ provided the initial energy is suitably small. However,    the boundary of $\Omega$  in \cite{Ho3} is flat.
Therefore, it remains completely open even for the existence of global weak solutions for any  $\ga\in (1,3/2]$ in general bounded domains. Hence,  it is interesting to study  the  global existence of weak and classical solutions to the initial-boundary-value problem   \eqref{a1}--\eqref{ch1} with any $\ga>1$ for general bounded smooth domains $\Omega $ with density containing vacuum initially.

Before stating the main results, we explain the notations and conventions used throughout this paper. We first give the definition  of simply connected domains.

\begin{definition}
Let $\Omega$ be a domain in $\r^3$. If the first Betti number of $\Omega$ vanishes, namely, any simple closed curve in $\Omega$ can be contracted to a point, we say that $\Omega$ is simply connected.  If the second Betti number of $\Omega$ is zero, we say that $\Omega$ has no holes.
\end{definition}

Then we define weak and classical solutions as follows.
\begin{definition}Let $T > 0$ be a finite constant.  A solution $(\rho,  u)$ to  \eqref{a1}  is called a weak solution if it satisfies \eqref{a1} in the sense of distribution. Moreover, when all the derivatives involved in \eqref{a1} are continuous functions, and \eqref{a1} holds  everywhere in $\Omega\times (0,  T)$, we call the solution a classical one.
\end{definition}

Next,   we set
$$ \int f dx\triangleq\int_{\Omega}fdx,$$
and
\bnn \bar{f}\triangleq\frac{1}{|\Omega|}\int_{\Omega}f  dx,\enn
which is the average of a function $f$ over $\Omega. $
   For   integer $k$ and $1\leq q<+\infty$, $W^{k,q}(\Omega)$ is the standard Sobolev spaces, $H^k(\Omega) \triangleq W^{k,2}(\Omega), $ and  $$ H_{  \o}^1(\Omega)\triangleq\left\{f\in
H^1 :  f\cdot n=0, ~~ {\rm curl}f\times n=-Af ~\mbox{ \rm on } ~\partial \Omega\right\}.$$
For some $s\in(0,1)$, the fractional Sobolev space $H^s(\Omega)$ is defined by
 $$ H^s(\Omega)\triangleq\left\{u\in L^2(\Omega)~\text{:} \int_{\Omega\times\Omega}\frac{|u(x)-u(y)|^2}{|x-y|^{n+2s}}dxdy<+\infty \right\},$$ which is a Banach space with the norm:
 $$\| u\|_{H^s(\Omega)}\triangleq \|u\|_{L^2(\Omega)}+\left(\int_{\Omega\times\Omega}\frac{|u(x)-u(y)|^2}{|x-y|^{n+2s}}dxdy\right)^\frac{1}{2}.$$

For simplicity, we denote $L^q(\Omega)$, $W^{k,q}(\Omega)$, $H^k(\Omega)$, and ${H^s(\Omega)}$ by $L^q$,  $W^{k,q}$, $H^k$,  and ${H^s}$ respectively.

For two $n\times n$  matrices $A=\{a_{ij}\},\,\,B=\{b_{ij}\}$, the symbol $A\colon  B$ represents the trace of $AB$, that is,
 $$ A\colon  B\triangleq \text{tr} (AB)=\sum\limits_{i,j=1}^{n}a_{ij}b_{ji}.$$


Finally, we   denote
the initial total energy of (\ref{a1})   as
\be \la{c0}
C_0 \triangleq\int_{\Omega}\left(\frac{1}{2}\n_0|u_0|^2 +G(\n_0) \right)dx,
\ee
with
\be\la{pn01}
G(\rho)\triangleq \rho \int_{\bar{\rho} }^{\rho}\frac{P(s)-P(\bar{\rho} )}{s^{2}} ds.
\ee
 Then one of the  main purposes of this paper is to establish the following global existence   of classical solutions of $\eqref{a1}$-$\eqref{ch1}$ in a general smooth bounded domain $\Omega\subset\r^3.$



\begin{theorem}\la{th1}  Let $\Omega$ be a simply connected bounded domain in $\r^3$ and its smooth boundary $\partial\Omega$ has a finite number of 2-dimensional connected components.  For given positive constants $M $ and   $\hat{\rho},$   suppose that the $3 \times 3$  symmetric matrix  $A $ in \eqref{ch1} is smooth and positive semi-definite, and the initial data $(\n_0,u_0)$ satisfy for some $q\in (3,6)$ and $s\in (1/2,1],$
\be \la{dt1}   (\rho_0 ,P(\rho_0) )  \in  W^{2,q}, \quad  u_0\in H^2\cap H^1_\o(\O) , \ee
\be\la{dt2} 0\leq\rho_0\leq\hat{\rho},~~\|u_0\|_{H^s}\leq M, \ee
and the compatibility condition
\be\la{dt3}
-\mu\triangle u_0-(\mu+\lambda)\nabla\div u_0 + \nabla P(\rho_0) = \rho_0^{1/2}g,
\ee
for some  $ g\in L^2.$
Then there exists a  positive constant  $\ve$   depending only on  $\mu$, $\lambda$, $\gamma$, $a$, $s$, $\hat{\rho}$, $M$, $\Omega,$ and the matrix $A$   such that if \be \la{c0q} C_0\le\ve\ee  with initial energy $C_0$ as in \eqref{c0},
the initial-boundary-value problem  \eqref{a1}-\eqref{ch1} has a unique classical solution $(\n,u)$ in $\Omega\times(0,\infty)$ satisfying that for any $0<\tau<T<\infty$,
\be\la{dt6}\begin{cases}
 (\rho ,P )\in C([0,T];W^{2,q} ),\\  \na u\in C([0,T];H^1 )\cap  L^\infty(\tau,T; W^{2,q}),\\
u_t\in L^{\infty} (\tau,T; H^2)\cap H^1 (\tau,T; H^1),\\   \sqrt{\n}u_t\in L^\infty(0,\infty;L^2) ,
\end{cases}\ee and that   for any $0< T<\infty$,
\be\la{dt5}
\ti C(T)\inf_{x\in \O}\n_0(x)\le\n(x,t)\le 2\hat{\n},\quad  (x,t)\in \O\times[0,T],
\ee
for some  positive constant $\ti C(T)$  depending only on $T$, $\mu$, $\lambda$, $\gamma$, $a$, $s$, $\hat{\rho}$, $M$, $\Omega,$ and the matrix $A.$
Moreover,  for any $r\in [1,\infty)$ and $p\in [1,6],$ there exist positive constants $C$ and $\eta_0$ depending only  on $\mu, \lambda, \gamma, a,s,\hat{\rho} , M, \bar\rho_0, \Omega, r, p,$  and the matrix $A$   such that for any $t\geq 1,$
\be  \la{qa1w}  \|\n-\bar\n_0\|_{L^r}+\|  u\|_{W^{1,p}} +\|\sqrt{\rho}\dot{u}\|^2_{L^2} \leq Ce^{-\eta_0 t}.\ee
\end{theorem}

Then, with the   exponential decay rate  \eqref{qa1w} at hand,  motivated by  the proof of \cite[Theorem 1.2]{lx},
we will   establish  the following large-time behavior of the spatial gradient of the
density when vacuum   appears initially.
\begin{theorem}\la{th2}
Under the conditions of Theorem \ref{th1}, assume further  that
there exists some point $x_0\in \Omega$ such that $\n_0(x_0)=0.$  Then the unique
global classical solution $(\n,u)$ to the   problem  \eqref{a1}-\eqref{ch1} obtained in
Theorem \ref{th1}  satisfies that for any $r_1>3,$   there exist positive constants $\hat{C}_1$ and $\hat{C}_2$ depending only  on $\mu$,  $\lambda$,  $\gamma$, $a$, $s$,   $\hat{\rho}$, $\bar\n_0,$ $M$, $\Omega$,  $r_1$ and the matrix $A$  such that for any $t\ge 1$,
\be\la{qa2w}\ba \|\na\n (\cdot,t)\|_{L^{r_1}}\geq \hat{C}_1 e^{\hat{C}_2 t} . \ea\ee
\end{theorem}

The third result concerns the global existence of weak solutions.
\begin{theorem}\label{th3} Under the conditions of Theorem \ref{th1} except \eqref{dt3}, where the condition \eqref{dt1} is replaced by
\be\la{dt7}\rho_0\in L^\infty(\Omega),\,\,  u_0\in H^1_\o(\O),\ee assume further that the initial energy $C_0$ as in \eqref{c0} satisfies \eqref{c0q} with $\ve$ as in Theorem \ref{th1}.
Then there exists at least one weak solution $(\rho,u)$ of the problem  \eqref{a1}-\eqref{ch1}   in
$\Omega\times(0,  \infty)$ satisfying \eqref{dt5},  \eqref{qa1w} and for any $0<\tau\leq T<\infty $ and $q\in [1,\infty)$,
\begin{equation}\label{111}
 \begin{cases}
\rho\in L^{\infty}(0, T; L^{\infty})\cap C([0,T];L^q),\\
u\in L^\infty(0, T;H^1), u_t\in L^2(\tau,T;L^6),\nabla u\in L^\infty(\tau,T;L^6),\\  {\rm curl}u,\,\, (2\mu+\lambda)\mathrm{div}u-P \in L^\infty(\tau,T;H^1  )\cap L^{2}(\tau,T;W^{1,6}).
       \end{cases}
\end{equation}

\end{theorem}

A few remarks are in order:

\begin{remark}\la{remark11} It should be mentioned here  that the    Navier-type slip condition  \eqref{Navi} is in fact  a particular case of   the slip boundary one \eqref{ch1}.
Indeed, since $u\cdot n=0$ on $\partial\Omega$, we have, for any unit tangential vector $\nu$,
\be\ba\la{ch7}
0=\frac{\partial}{\partial\nu}(u\cdot n)&=(D(u)\,n)\cdot\nu-\frac{1}{2}\curl u\times n\cdot\nu + \nu\cdot\nabla n\cdot u\\
&=(D(u)\,n)\cdot\nu-\frac{1}{2}\curl u\times n\cdot\nu + u\cdot D(n)\cdot \nu,\ea\ee
where in the last equality we have used the fact $\curl n\cdot n=0$. Consequently, \eqref{Navi} is equivalent to
\be \la{ch5}
\curl u\times n=(2D(n)-\vartheta I)u\ee
where $I$ is $3\times 3$ identity matrix and  $\vartheta I-2D(n)$ is indeed a $3\times 3$ symmetric matrix.
\end{remark}

\begin{remark} Since $q>3,$ it follows from Sobolev's inequality and \eqref{dt6}$_1$  that \be\la{soh1}  \n,\na \n \in C(\bar\Omega\times [0,T]).\ee
Moreover, it also follows from \eqref{dt6}$_2$ and \eqref{dt6}$_3$ that  \be \la{soh2} u,\na u, \na^2 u, u_t \in C(\bar\Omega\times [\tau,T]),\ee due to the following simple fact that $$L^2(\tau,T;H^1)\cap H^1(\tau,T;H^{-1})\hookrightarrow C([\tau,T];L^2).$$
Finally, by \eqref{a1}$_1,$ we have \bnn \n_t=-u\cdot\na \n-\n\div u\in C(\bar\Omega\times [\tau,T]),\enn which together with \eqref{soh1} and \eqref{soh2} shows that the solution obtained by Theorem \ref{th1} is a classical one.
\end{remark}
\begin{remark}
Theorem \ref{th2} implies that   the  oscillation of the density will grow unboundedly in the long run with an exponential 
  rate provided  vacuum (even a point) appears initially.  This new phenomena is   somewhat surprisingly compared with the  Cauchy problem  (\cite{hlx1,jx01}) where  there is   not any result concerning the growth rate of   the  gradient of the density.
\end{remark}

\begin{remark}\la{remark4} For the sake of simplicity, we assume that the matrix $A$ is smooth and  positive semi-definite. However, these conditions can be relaxed. Indeed, we only use the assumption that  the matrix $A$ is   positive semi-definite in the proof of \eqref{a16} and \eqref{ax3999}(see \eqref{m8} and \eqref{ax4000}). Thus, let  $\lambda_i(x) (i=1,2,3)$ be the eigenvalues of $A$  whose negative parts are denoted by $\lambda_1^-(x)$, $\lambda_2^-(x)$, and $\lambda_3^-(x)$ respectively. Then  one can deduce that  \eqref{a16} and \eqref{ax3999} both still hold provided  $\lambda_1^-(x), \lambda_2^-(x), \lambda_3^-(x)$ are bounded by some suitably small positive constant   depending  only  on $\lambda$, $\mu$, the constants of Poincar\'{e}'s inequality and the constant $C_1$ in \eqref{tdu1} for $p=2$. The other restriction on $A$ comes from a priori estimates related to A, in which \eqref{remark2} plays a decisive role. In fact, by Lemma \ref{crle1}, \eqref{x2666}, \eqref{x268} and \eqref{h18}, one can find that for any $p\in [2,6]$,
 \be\la{remark3}\ba
\|\nabla^{3}u\|_{L^p}&\leq C(\|\div u\|_{W^{2,p}}+\|\curl u\|_{W^{2,p}})\\
&\leq C(\|\rho\dot{u}\|_{W^{1,p}}+\|\nabla(Au)^\perp\|_{W^{1,p}}+\|\nabla u\|_{L^2})\\&\quad+C\|\nabla P\|_{W^{1,p}}+C\|P-\Bar P\|_{L^p} .
\ea\ee
 Therefore, compared with \eqref{remark2}, it is sufficient to assume that $A\in W^{2,6}$ rather than smooth. Moreover, for the lower order priori estimates in Section \ref{se3}, it is enough to suppose that $A\in W^{1,6}$, and then the dependence of $\ve$ on $A$ in Theorem \ref{th1} can be determined by $\|A\|_{W^{1,6}}$.
\end{remark}

The following  result   removes   the condition that the region is simply connected  provided $2\mu +  3 \lambda>0.$ 

\begin{theorem}\label{th4}  Assume that $\Omega$ is a smooth bounded  domain in $\r^3$, not necessarily simply connected, and its smooth boundary $\partial\Omega$ has a finite number of 2-dimensional connected components. Let $A\in W^{2,6}(\Omega)$ satisfy  that $A+2D(n)$  is a positive semi-definite $3\times 3$ symmetric matrix. Moreover,  $A+2D(n) $   is  positive on some $\Sigma\subset\partial\Omega$ with $|\Sigma|>0$  when $\Omega$ is axially symmetric. Then, the conclusions of Theorems \ref{th1}-\ref{th3} are still valid provided $2\mu +  3 \lambda>0.$
\end{theorem}

By virtue of \eqref{ch5}, $A=\vartheta I-2D(n)$ corresponds to Navier-type slip boundary condition \eqref{Navi}.    As a direct consequence of  Theorems \ref{th1}-\ref{th4} (see their remarks also),  for compressible Navier-Stokes equations \eqref{a1} with Navier-type slip boundary condition \eqref{Navi}, we have the following conclusion on the global existence and large-time behavior of classical or weak solutions.
\begin{corollary} Let $\Omega$ be a smooth bounded  domain in $\r^3$ and its smooth boundary $\partial\Omega$ has a finite number of 2-dimensional connected components.
Then, for $\vartheta\in W^{2,6}(\Omega),$  
 the conclusions of Theorems \ref{th1}-\ref{th3} still hold where  \eqref{ch1} is replaced by  \eqref{Navi} provided    one of the following two conditions holds:

(1)  $\Omega$ is simply connected, $A=\vartheta I-2D(n)$ satisfies the assumption given by Remark \ref{remark4};

(2)  $2\mu+3\lambda> 0$, $\vartheta\geq 0$ and in addition, $\vartheta>0$ on some $\Sigma\subset\partial\Omega$ with $|\Sigma|>0$ if $\Omega$ is axially symmetric.
\end{corollary}
\begin{remark} Compared with \cite[Theorem 7.69]{ANIS} where they obtain the global weak solutions to the compressible Navier-Stokes equations with Navier-slip boundary conditions ($\vartheta=0$) for $\ga>3/2$ in a non axisymmetric domain,  our results    establish  the global existence of both weak and classical solutions (which may contain vacuum states) for any $\ga>1$ and more relaxed conditions for $\vartheta$ provided that the initial energy is suitably small. Moreover, \eqref{dt5} indicates  that vacuum states  will not exhibit in any finite time even for weak solutions provided that no
vacuum states are present initially.\end{remark}

We now comment on the analysis of this paper. Indeed, compared with the previous results (\cite{hlx1,jx01}) where they treated the Cauchy problem,  the slip boundary condition \eqref{ch1} causes additional difficulties in developing a priori estimates for   solutions of the compressible Navier-Stokes equations. To overcome the difficulties,  our research bases on three observations.    First, it is important to find an equivalence of norms in $H^1$. Thanks to \cite{vww}, for simply connect bounded domains, we have
$$\|\nabla u\|_{L^2}\leq C(\|\div u\|_{L^2}+\|\curl u\|_{L^2}).$$
And for more general bounded domains, we have \eqref{equiv1}. Next, for $v=(v^1,v^2,v^3),$ denoting the material derivative of $v$ by \be\la{oio1}\dot v\triangleq v_t+u\cdot\nabla v,\ee    we rewrite $ (\ref{a1})_2 $ in the form
\be \la{hod1}\ba
\rho\dot{u}=\nabla F - \mu\nabla\times\curl u ,
\ea \ee with
\be \la{dt0}  \text{curl} u \triangleq \nabla\times u ,\quad F\triangleq(\lambda+2\mu)\,\div u-(P-\Bar P),\ee where  the vorticity $\text{curl} u $   and the so-called the effective viscous flux $F$     both play  an important role in our following analysis.  Since $u\cdot n=0$  on $\partial\Omega$, we check that
\be\la{pzw1} u\cdot\nabla u\cdot n=-u\cdot\nabla n\cdot u,\ee which implies (see \eqref{paz2})
\bnn(\dot{u}+(u\cdot\nabla n)\times u^{\perp})\cdot n=0 \mbox{ on } \pa\O  ,\enn  with $u^{\perp}\triangleq -u\times n$ on $\partial\Omega.$ As a direct consequence of this observation, we have (see \eqref{tb90})
\bnn
\ba\|\dot{u}\|_{L^6}\le C(\|\nabla\dot{u}\|_{L^2}+\|\nabla u\|_{L^2}^2).
\ea\enn
Similarly, one can get $(\curl u+(Au)^\perp)\times n=0$ on $\partial\Omega$ by the other boundary condition $\curl u\times n=-A\, u$.  Combining this with \eqref{hod1}   implies that one can  treat $(\ref{a1})_2$ as a Helmholtz-Wyle decomposition of $\rho\dot{u}$ which   makes it possible to estimate $\nabla F$ and $\nabla \curl u$ (see \eqref{x266}). Finally, since $u\cdot n=0$  on $\partial\Omega$, we have  \bnn  u=   u^\bot\times n \mbox{ on } \pa\O,\enn which, combined with   the   simple fact that  $\div(   \na u^i \times  u^\bot)=-\na u^i\cdot \na\times u^\bot, $ implies that we can bound the following key boundary integrals  concerning the effective viscous flux $F,$ the vorticity ${\rm curl}u,$ and $\na u$(see \eqref{bz3} and \eqref{ax4000}) $$\int_{\partial\Omega}F  ({u}\cdot\na)u\cdot\na n\cdot uds, \,\,
 \int_{\partial\Omega}\curl u_t\times n\cdot\dot{u}ds.$$   All these treatments are the key to estimating the crucial integrals on the boundary $\partial\Omega$.

The rest of the paper is organized as follows. First,  some notations, known facts and elementary inequalities needed in later analysis are collected in Section 2. Sections 3 and  4 are devoted to deriving the necessary a priori estimates on classical solutions which can guarantee the extension of the local classical solution to be a global one. Finally, the main results, Theorems \ref{th1}--\ref{th4} will be proved in Sections 5 and 6.

\section{Preliminaries}\la{se2}

In this section, we recall some  known facts and elementary inequalities which will be used later.


First, similar to the proof of \cite[Theorem 1.4]{hxd1}, we have the local existence of strong and classical solutions.
\begin{lemma}\la{loc1} Let $\Omega$ be as in Theorem \ref{th1}, assume that $(\n_0,u_0)$ satisfies \eqref{dt1} and \eqref{dt3}. Then there exist a small time $T>0$ and a unique strong solution $(\n,u)$ to the problem \eqref{a1}-\eqref{ch1} on $\Omega\times(0,T]$ satisfying for any $ \tau\in(0,T),$
\be\nonumber\begin{cases}
 (\rho ,P )\in C([0,T];W^{2,q} ),\\  \na u\in C([0,T];H^1 )\cap  L^\infty(\tau,T; W^{2,q}),\\
u_t\in L^{\infty} (\tau,T; H^2)\cap H^1 (\tau,T; H^1),\\   \sqrt{\n}u_t\in L^\infty(0,T;L^2) .
\end{cases}\ee \end{lemma}


 Next,  the well-known Gagliardo-Nirenberg's inequality (see \cite{nir})
  will be used frequently later.
\begin{lemma}
[Gagliardo-Nirenberg]\la{l1} Assume that $\Omega$ is a bounded Lipschitz domain in $\r^3$. For  $p\in [2,6],\,q\in(1,\infty), $ and
$ r\in  (3,\infty),$ there exist    generic
 constants
$C, C_1,C_2 $ which   depend only on $p$, $q$, $r$, and $\Omega$ such that for any  $f\in H^1({\O }) $
and $g\in  L^q(\O )\cap D^{1,r}(\O), $
\be\la{g1}\|f\|_{L^p(\O)}\le C  \|f\|_{L^2}^{\frac{6-p}{2p}}\|\na
f\|_{L^2}^{\frac{3p-6}{2p}}+C_1\|f\|_{L^2} ,\ee
\be\la{g2}\|g\|_{C\left(\ol{\O }\right)} \le C
\|g\|_{L^q}^{q(r-3)/(3r+q(r-3))}\|\na g\|_{L^r}^{3r/(3r+q(r-3))} + C_2\|g\|_{L^2}.
\ee
Moreover, if either $f\cdot n|_{\partial\Omega}=0 $ or $\bar{f}=0,$     we can choose $C_1=0.$ Similarly,   the constant $C_2=0 $ provided $ g\cdot n|_{\partial\Omega}=0$ or $\bar{g}=0$.
\end{lemma}

In order to get the
uniform (in time) upper bound of the density $\n,$ we need the following Zlotnik's  inequality.
\begin{lemma}[\cite{zl1}]\la{le1}   Suppose the function $y$ satisfies
\bnn y'(t)= g(y)+b'(t) \mbox{  on  } [0,T] ,\quad y(0)=y^0, \enn
with $ g\in C(R)$ and $y, b\in W^{1,1}(0,T).$ If $g(\infty)=-\infty$
and \be\la{a100} b(t_2) -b(t_1) \le N_0 +N_1(t_2-t_1)\ee for all
$0\le t_1<t_2\le T$
  with some $N_0\ge 0$ and $N_1\ge 0,$ then
\bnn y(t)\le \max\left\{y^0,\zeta_0 \right\}+N_0<\infty
\mbox{ on
 } [0,T],
\enn
  where $\zeta_0 $ is a constant such
that \be\la{a101} g(\zeta)\le -N_1 \quad\mbox{ for }\quad \zeta\ge \zeta_0.\ee
\end{lemma}

Next, consider the Lam\'{e}'s system
\be\la{cxtj1}\begin{cases}
-\mu\Delta u-(\lambda+\mu)\nabla\div u=f, \,\, &x\in\Omega, \\
u\cdot n=0,\,\curl u\times n=-Au,\,\,&x\in\partial\Omega,
\end{cases} \ee
where $u=(u^{1},u^{2},u^{3}),\,\,f=(f^{1},f^{2},f^{3})$, $\Omega$ is a bounded smooth domain in $\r^3,$ and $\mu,\lambda$ satisfy the condition \eqref{h3}.
It follows from \cite{soln2} that the Lam\'{e}'s system is of Petrovsky type. In Petrovsky's systems, roughly speaking, different equations and unknowns have the same ``differentiability order", see \cite{soln1}. We also recall that Petrovsky's systems belong to an important subclass of Agmon-Douglis-Nirenberg (ADN) elliptic systems (see\cite{adn}), which has the same good properties of self-adjoint ADN systems. 
Thus,   we have  the following  standard estimates.
\begin{lemma}  [\cite{adn}] \la{zhle}
Let $u$ be a smooth solution of the Lam\'{e}'s equation \eqref{cxtj1}. Then for   $q\in(1,\infty),\,\, k\geq0,$ there exists a positive constant $C$ depending only on $\lambda,\,\mu,\,q,\,\,k$, $\Omega$ and the matrix $A$ such that

(1) If $f\in W^{k,q},$   then
$$\|u\|_{W^{k+2,q}}\leq C(\|f\|_{W^{k,q}}+\|u\|_{L^q}),$$

(2) If $f=\nabla g$ and $g\in W^{k,q},$  then
$$\|u\|_{W^{k+1,q}}\leq C(\|g\|_{W^{k,q}}+\|u\|_{L^q}).$$
\end{lemma}

Next, the following two lemmas can be found in   \cite[Theorem 3.2]{vww} and   \cite[Propositions 2.6-2.9]{CANEHS}.
\begin{lemma} \la{crle1}
Let $k\geq0$ be a integer and  $\Omega$ be a bounded domain in $\r^3$ with $C^{k+1,1}$ boundary $\partial\Omega$, $1<q<+\infty$. Then for $v\in W^{k+1,q}$ with $v\cdot n=0$ on $\partial\Omega$, there exists a constant $C=C(q,k,\Omega)$ such that

\be\la{paz11} \|v\|_{W^{k+1,q}}\leq C(\|\div v\|_{W^{k,q}}+\|\curl v\|_{W^{k,q}}+\|v\|_{L^q}).\ee

Furthermore, assume that $\Omega$ is a simply connected bounded domain in $\r^3$ with $C^{k+1,1}$ boundary $\partial\Omega$, we have
\be \la{paz1}\|v\|_{W^{k+1,q}}\leq C(\|\div v\|_{W^{k,q}}+\|\curl v\|_{W^{k,q}}).\ee
\end{lemma}
\begin{lemma}   \la{crle2}
Let $k\geq0$ be an  integer, $1<q<+\infty$. Suppose that $\Omega$ is a bounded domain in $\r^3$ and its $C^{k+1,1}$ boundary $\partial\Omega$ only has a finite number of 2-dimensional connected components. Then for $v\in W^{k+1,q}$ with $v\times n=0$ on $\partial\Omega$, there exists a constant $C=C(q,k,\Omega)$ such that
$$\|v\|_{W^{k+1,q}}\leq C(\|\div v\|_{W^{k,q}}+\|\curl v\|_{W^{k,q}}+\|v\|_{L^q}).$$
In particular, if  $\Omega$ has no holes, then
$$\|v\|_{W^{k+1,q}}\leq C(\|\div v\|_{W^{k,q}}+\|\curl v\|_{W^{k,q}}).$$
\end{lemma}

Next,  to estimate the $L^1(0,T;L^\infty(\O))$-norm of  $\nabla u  ,$  
 we need the following    Beale-Kato-Majda type inequality with respect to the slip boundary condition \eqref{ch1}, which was first proved in \cite{bkm,kato} in the whole 3D spatial space when $\div u\equiv 0.$
\begin{lemma}\la{le9} Let $\Omega$ be a bounded domain in $\r^3$ with smooth boundary.
For $3<q<\infty$, assume that $u\cdot n=0$ and $\curl u\times n=-Au$ on $\partial\Omega$, $ u\in W^{2,q}$, then there is a constant  $C=C(q, \Omega, A)$ such that  the following estimate holds
\be\la{zhbc0}\ba
\|\na u\|_{L^\infty}\le C\left(\|{\rm div}u\|_{L^\infty}+\|\curl u\|_{L^\infty} \right)\ln(e+\|\na^2u\|_{L^q})+C\|\na u\|_{L^2} +C .
\ea\ee
\end{lemma}
\begin{proof} We borrow some ideas   of \cite{hlx,hl01} and make some slight modifications.
It follows from \cite{soln1,soln2} that $u$ can be represented in the form
\bnn\ba u^{i}&=\int G_{i,\cdot}(x,y)\cdot(\mu\Delta_y u+(\lambda+\mu)\nabla_y\div_y u)dy\\
 &\triangleq\int G_{ij}(x,y)(\mu\Delta_y u^j+(\lambda+\mu)\partial_{y^j}\div_y u)(y)dy,\ea\enn
where, $G=\{G_{ij}\}$ with $G_{ij}=G_{ij}(x,y)\in C^{\infty}(\Omega\times\Omega\backslash D)$, $D\equiv\{(x,y)\in\Omega\times\Omega:\,x=y\}$, is Green matrix of the Lam\'{e}'s system
\eqref{cxtj1} and satisfies that for every multi-indexes $\alpha=(\alpha^1,\alpha^2,\alpha^3)$ and $\beta=(\beta^1,\beta^2,\beta^3)$, there is a constant $C_{\alpha,\beta}$ such that for all
$(x,y)\in\Omega\times\Omega \backslash D,$ and $i,j=1,2,3,$
$$|\partial_{x}^{\alpha}\partial_{y}^{\beta}G_{ij}(x,y)|\leq C_{\alpha,\beta}|x-y|^{-1-|\alpha|-|\beta|},$$
here $|\alpha|=\alpha^1+\alpha^2+\alpha^3$ and $|\beta|=\beta^1+\beta^2+\beta^3$.

Notice that according to the definition of $A$ in \eqref{th1}, $A u$ is still a tangential vector on $\partial\Omega$, and then we set \be\la{ljq1}(Au)^\perp\triangleq-(Au)\times n, \ee so $Au=(Au)^\perp\times n$.
Therefore,
\be\la{zhbc1}\ba
u^{i}(x)=&(\lambda+2\mu)\int G_{i,\cdot}(x,y)\cdot\nabla_y\div u(y)dy\\&-\mu\int G_{i,\cdot}(x,y)\cdot\nabla_y\times\curl u(y)dy \\
=&(\lambda+2\mu)\int G_{i,\cdot}(x,y)\cdot\nabla_y\div u(y)dy\\&-\mu\int G_{i,\cdot}(x,y)\cdot\nabla_y\times(\curl u+(Au)^\perp)dy \\
&  +\mu\int\nabla_y\times (A u(y))^\perp\cdot G_{i,\cdot}(x,y) d y \triangleq\sum_{j=1}^{3}U_j^i.
\ea\ee
It suffices to estimate the three terms $U_j^i,\,j=1,2,3$. Let $\delta\in(0,1]$ be a constant to be chosen and introduce a cut-off function $\eta_{\delta}(x)$ satisfying $\eta_{\delta}(x)=1$ for $|x|<\delta;$ $\eta_{\delta}(x)=0$ for $|x|>2\delta,$ and $|\nabla\eta_{\delta}(x)|<C\delta^{-1}.$ Notice that $G_{i,\cdot}(x,y)\cdot n=0$ on $\partial\Omega$, $\nabla U_1^i$ can be written as
\be\la{zhbc2}\ba
\nabla U_1^i&=(\lambda+2\mu)\int\eta_{\delta}(|x-y|)\,\nabla_x G_{i,\cdot}(x,y)\nabla_y\div u(y)dy\\
&\quad+(\lambda+2\mu)\int\nabla_y\eta_{\delta}(|x-y|)\cdot\nabla_x G_{i,\cdot}(x,y)\,\div u(y)dy \\
&\quad-(\lambda+2\mu)\int(1-\eta_{\delta}(|x-y|))\,\nabla_x\div_y G_{i,\cdot}(x,y)\,\div u(y)dy \\
&\triangleq (\lambda+2\mu)\sum_{k=1}^{3}\tilde{I}_k.
\ea\ee
Now we estimate $\tilde{I}_k,\,\,k=1,2,3.$
\be\la{zhbc3}\ba
|\tilde{I}_1|&\leq C\|\eta_{\delta}(|x-y|)\,\nabla_x G_{i,\cdot}(x,y)\|_{L^{q/(q-1)}}\|\nabla^{2}u\|_{L^q}\\
&\leq C\left(\int_0^{2\delta}r^{-2q/(q-1)}r^{2}dr\right)^{(q-1)/q}\|\nabla^{2}u\|_{L^q}\\
&\leq C\delta^{(q-3)/q}\|\nabla^{2}u\|_{L^q},
\ea\ee
\be\la{zhbc4}\ba
|\tilde{I}_2|&=\left|\int\nabla_y\eta_{\delta}(|x-y|)\cdot\nabla_x G_{i,\cdot}(x,y)\,\div u(y)dy\right|\\
&\leq C\int|\nabla_y\eta_{\delta}(y)\cdot\nabla_x  G_{i,\cdot}(x,y)|dy\,\|\div u\|_{L^{\infty}}\\
&\leq C\int_\delta^{2\delta}\delta^{-1}r^{-2}r^{2}dr\,\|\div u\|_{L^{\infty}}\\
&\leq C \|\div u\|_{L^{\infty}},
\ea\ee
\be\la{zhbc5}\ba
|\tilde{I}_3|&=\left|\int(1-\eta_{\delta}(|x-y|))\,\nabla_x\div_y G_{i,\cdot}(x,y)\,\div u(y)dy\right|\\
&\leq C\left(\int_{\delta\leq|x-y|\leq1}+\int_{|x-y|>1}\right)|\nabla_x\div_y G_{i,\cdot}(x,y)|\,|\div u(y)|dy\\
&\leq C\int_\delta^{1}r^{-3}r^{2}dr\,\|\div u\|_{L^{\infty}}+C\left(\int_1^{\infty}r^{-6}r^2dr\right)^{\frac{1}{2}}\|\div u\|_{L^2}\\
&\leq -C\ln\delta\,\|\div u\|_{L^{\infty}}+C\|\nabla u\|_{L^2}.
\ea\ee
It follows from \eqref{zhbc2}-\eqref{zhbc5} that
\be\la{zhbc6}\ba
\|\nabla U_1\|_{L^{\infty}}\leq C\left(\delta^{(q-3)/q}\|\nabla^{2}u\|_{L^q}+(1-\ln\delta)\,\|\div u\|_{L^{\infty}}+\|\nabla u\|_{L^2}\right).
\ea\ee
Since by \eqref{ch1}, $(\curl u+(Au)^\perp)\times n=0$ on $\partial\Omega,$   we  rewrite $\nabla U_2^i$ as
\bnn\ba
\nabla U_2^i&=-\mu\int\eta_{\delta}(|x-y|)\,\nabla_x G_{i,\cdot}(x,y)\cdot\nabla_y\times(\curl u+(Au)^\perp)dy\\
&\quad+\mu\int\nabla_y\eta_{\delta}(|x-y|)\times\nabla_x G_{i,\cdot}(x,y)\cdot(\curl u+(Au)^\perp)dy \\
&\quad-\mu\int(1-\eta_{\delta}(|x-y|))\,\nabla_y\times\nabla_x G_{i,\cdot}(x,y)\,\cdot(\curl u+(Au)^\perp)dy. \\
\ea\enn
A discussion similar to the previous term gives
\be\la{zhbc7}\ba
 \|\nabla U_2\|_{L^{\infty}} \leq& C \delta^{(q-3)/q}\|\nabla^{2}u\|_{L^q}\\&+C(1-\ln\delta)\,(\|\curl u\|_{L^{\infty}}+\|u\|_{L^{\infty}})+C\|\nabla u\|_{L^2} .
\ea\ee

Finally, it is clear that
\be\la{zhbc8}\ba
\|\nabla U_3\|_{L^{\infty}}\leq C\left(\delta^{(q-3)/q}\|\nabla^{2}u\|_{L^q}+(1-\ln\delta)\,\|u\|_{L^{\infty}}+\|\nabla u\|_{L^2}\right).
\ea\ee
Combining \eqref{zhbc1} with \eqref{zhbc6}-\eqref{zhbc8}, and utilizing \eqref{paz11} leads to
\bnn\ba
\|\nabla u\|_{L^{\infty}}\leq& C \delta^{(q-3)/q}\|\nabla^{2}u\|_{L^q}+C(1-\ln\delta) (\|\curl u\|_{L^{\infty}}+\|\div u\|_{L^{\infty}})\\&+C\|\nabla u\|_{L^2} ,
\ea\enn
which after choosing  $\delta=\min\{1,\|\nabla^{2}u\|_{L^q}^{-q/(q-3)}\}$   gives \eqref{zhbc0} and completes  the proof of Lemma \ref{le9}.
\end{proof}

Next,
for the problem
\be\la{e480}\begin{cases}
{\rm div}v=f,&x\in\Omega, \\
v=0,&x\in{\partial\Omega},
\end{cases} \ee one has the following conclusion.
\begin{lemma} \cite[Theorem III.3.1]{GPG} \la{th00}  There exists a linear operator $\mathcal{B} = [\mathcal{B}_1 , \mathcal{B}_2 , \mathcal{B}_3 ]$ enjoying
the properties:

1)The operator $$\mathcal{B}:\{f\in L^p(\O):\bar f =0\}\mapsto (W^{1,p}_0(\O))^3$$ is a bounded linear one, that is,
\bnn \|\mathcal{B}[f]\|_{W^{1,p}_0(\O)}\le C(p)\|f\|_{L^p(\O)}, \mbox{ for any }p\in (1,\infty).\enn

2) The function $v = \mathcal{B}[f]$ solves the problem \eqref{e480}.

3) If, moreover, for $f = \div  g$ with a certain $g\in L^r(\O),$ $g\cdot n|_{\pa\O}=0,$  then  for   any $r  \in (1,\infty),$
\bnn \|\mathcal{B}[f]\|_{L^{r}(\O)}\le C(r)\|g\|_{L^r(\O)} .\enn
\end{lemma}

Now, for $F$, $\curl u$ as in \eqref{dt0}, we have the following    key   a priori estimates which will be used frequently.
\begin{lemma}   \la{le3}
 Assume $\Omega$ is a bounded domain in $\r^3$ and its smooth boundary $\partial\Omega$ only has a finite number of 2-dimensional connected components. Let $(\rho,u)$ be a smooth solution of \eqref{a1} with slip boundary condition \eqref{ch1}. Then for any $p\in[2,6]$, there exists a positive constant $C$ depending only on  $p$, $\mu$, $\lambda$ and $A$ such that
\be\la{h20}\ba
\|F\|_{L^p}&\leq C\|\rho\dot{u}\|_{L^2}^{(3p-6)/(2p)}(\|\nabla u\|_{L^2}+\|P-\Bar P\|_{L^2})^{(6-p)/(2p)}\\
&\quad+C(\|\nabla u\|_{L^2}+\|P-\Bar P\|_{L^2}),
\ea\ee
\be\la{hh21}\ba
\|\curl u\|_{L^p}\leq C\|\rho\dot{u}\|_{L^2}^{(3p-6)/(2p)}\|\nabla u\|_{L^2}^{(6-p)/(2p)}+C\|\nabla u\|_{L^2},
\ea\ee
 \be\la{hh22}\ba
\|F\|_{L^p}+\|\curl u\|_{L^p}\leq C(\|\rho\dot{u}\|_{L^2}+\|\nabla u\|_{L^2}),
\ea\ee
\be\la{h18}\ba
\|\nabla u\|_{L^p}&\leq C\|\rho\dot{u}\|_{L^2}^{(3p-6)/(2p)}(\|\nabla u\|_{L^2}+\|P-\Bar P\|_{L^2})^{(6-p)/(2p)}\\
&\quad+C(\|\nabla u\|_{L^2}+\|P-\Bar P\|_{L^p}),
\ea\ee
\be\la{h19}\ba
\|\nabla F\|_{L^p}+\|\nabla\curl u\|_{L^p}\leq C(\|\rho\dot{u}\|_{L^p}+\|\nabla u\|_{L^2}+\|P -\Bar P   \|_{L^p}).
\ea\ee
\end{lemma}
\begin{proof} For $(Au)^\perp $ as in \eqref{ljq1} and $(\curl u+(Au)^\perp)\times n=0$ on $\partial\Omega$, we have, for any $\eta\in C^{\infty}(\r^3)$,
\bnn\la{curl1}\ba&
 \int\nabla\times\curl u \cdot\nabla\eta dx\\&=\int\nabla\times(\curl u+(Au)^\perp)\cdot\nabla\eta dx-\int\nabla\times(Au)^\perp \cdot\nabla\eta dx\\
&=-\int\nabla\times(Au)^\perp \cdot\nabla\eta dx,
\ea\enn which together with   $(\ref{a1})_2 $ yields that the viscous flux $F$ satisfies that for any $\eta\in C^{\infty}(\r^3)$,
$$\int\nabla F\cdot\nabla\eta dx=\int\left(\rho\dot{u}-\nabla\times(Au)^\perp\right)\cdot\nabla\eta dx, $$
that is,
\bnn\begin{cases}
\Delta F=\div(\rho\dot{u}-\nabla\times(Au)^\perp),~~ &x\in\Omega,\\ \frac{\partial F}{\partial n}=(\rho\dot{u}-\nabla\times(Au)^\perp)\cdot n,\,\, &x\in\partial\Omega.
\end{cases}\enn
It follows from \cite[Lemma 4.27]{ANIS} that
\be\la{x266}\ba
\|\nabla F\|_{L^q}&\leq C(\|\rho\dot{u}\|_{L^q}+\|\nabla\times(Au)^\perp\|_{L^q})\\
&\leq C(\|\rho\dot{u}\|_{L^q}+\|\nabla u\|_{L^q}),
\ea\ee
and that for any integer $k\geq 0$,
\be\la{x2666}\ba
\|\nabla F\|_{W^{k+1,q}}\leq C(\|\rho\dot{u}\|_{W^{k+1,q}}+\|\nabla\times(Au)^\perp\|_{W^{k+1,q}}).
\ea\ee

Furthermore, since $\bar F=0$, one  deduces from \eqref{g1} and \eqref{x266} that for $p\in[2,6]$,
\be\la{ui1}\ba
\|F\|_{L^p}\leq C\|\nabla F\|_{L^2}\leq C(\|\rho\dot{u}\|_{L^2}+\|\nabla u\|_{L^2}),
\ea\ee and
\bnn\ba
\|F\|_{L^p}&\leq C\|F\|_{L^2}^{(6-p)/(2p)}\|\nabla F\|_{L^2}^{(3p-6)/(2p)}\\
&\leq C\|\rho\dot{u}\|_{L^2}^{(3p-6)/(2p)}(\|\nabla u\|_{L^2}+\|P-\Bar P\|_{L^2})^{(6-p)/(2p)}\\
&\quad+C\|\nabla u\|_{L^2}+C\|P-\Bar P\|_{L^2},
\ea\enn which is \eqref{h20}.

One rewrites $(\ref{a1})_2 $ as $\mu\nabla\times\curl u=\nabla F-\rho\dot{u}.$
Since $(\curl u+(Au)^\perp)\times n=0$ on $\partial\Omega$ and $\div(\nabla\times\curl u)=0$, by Lemma \ref{crle2}, we get
\be\la{x267}\ba
\|\nabla\curl u\|_{L^q}&\leq C(\|\nabla\times\curl u\|_{L^q}+\|\nabla u\|_{L^q})\\
&\leq C(\|\rho\dot{u}\|_{L^q}+\|\nabla u\|_{L^q}),
\ea\ee
and for any integer $k\geq0$,
\be\la{x268}\ba
\|\nabla\curl u\|_{W^{k+1,q}}&\leq C(\|\nabla\times\curl u\|_{W^{k+1,q}}+\|\curl u\|_{L^q}+\|(Au)^\perp\|_{W^{k+2,q}})\\
&\leq C(\|\rho\dot{u}\|_{W^{k+1,q}}+\|\nabla(Au)^\perp\|_{W^{k+1,q}}+\|\nabla u\|_{L^q}),
\ea\ee
where we have taken advantage of \eqref{x2666}. Thus,  by Gagliardo-Nirenberg's inequality and \eqref{x267}, we get \eqref{hh21}. Combining   \eqref{hh21} with \eqref{ui1}  shows \eqref{hh22}.

Next, by virtue of \eqref{paz11}, \eqref{h20} and \eqref{hh21}, it indicates that
\bnn\ba
\|\nabla u\|_{L^p}&\le C(\|\div u\|_{L^p}+\|\curl u\|_{L^p}+ \|u\|_{L^p}) \\
&\le C(\|F\|_{L^p}+\|\curl u\|_{L^p}+\|P-\Bar P\|_{L^p}+\|\nabla u\|_{L^2} ) \\
&\le C\|\rho\dot{u}\|_{L^2}^{(3p-6)/(2p)}(\|\nabla u\|_{L^2}+\|P-\Bar P\|_{L^2})^{(6-p)/(2p)} \\
&\quad+C(\|\nabla u\|_{L^2}+\|P-\Bar P\|_{L^p}).
\ea\enn
which together with \eqref{x266} and \eqref{x267} gives \eqref{h18}.

The inequality \eqref{h19} is a direct consequence of \eqref{x266} and  \eqref{x267} and  we complete the proof of Lemma \ref{le3}.
\end{proof}
Finally, using the boundary condition \eqref{ch1}, we have the following estimates on the material derivative of $u.$

\begin{lemma}\la{uup1} Under the assumption of Lemma \ref{le3}, there exists a positive constant $\Lambda$ depending only on $\Omega$ such that
\be\la{tb90}
\ba\|\dot{u}\|_{L^6}\le \Lambda(\|\nabla\dot{u}\|_{L^2}+\|\nabla u\|_{L^2}^2),
\ea\ee
\be\la{tb11}\ba
\|\nabla\dot{u}\|_{L^2}\le \Lambda(\|\div \dot{u}\|_{L^2}+\|\curl \dot{u}\|_{L^2}+\|\nabla u\|_{L^4}^2).
\ea\ee
\end{lemma}

\begin{proof} First, setting $u^{\perp}\triangleq -u\times n,$  we have by \eqref{ch1}
\be\la{ie01}\dot{u}\cdot n= u\cdot\nabla u\cdot n=-u\cdot\nabla n\cdot u=-(u\cdot\nabla n)\times u^{\perp}\cdot n\,\,\,\text{on}\,\,
\partial\Omega,\ee due to the following simple fact $$  v\times (u\times n  )=(v\cdot n)u-(v\cdot u)n,$$ with $v=u\cdot \na n.$ It thus follows for \eqref{ie01} that
\be\la{paz2}(\dot{u}+(u\cdot\nabla n)\times u^{\perp})\cdot n=0 \mbox{ on } \partial\Omega,\ee which together with Poincar\'{e}'s inequality gives
$$\|\dot{u}+(u\cdot\nabla n)\times u^{\perp}\|_{L^\frac{3}{2}}\le C\|\nabla(\dot{u}+(u\cdot\nabla n)\times u^{\perp} )\|_{L^\frac{3}{2}}.$$ Thus, we obtain
\bnn\la{tb9}
\ba\|\dot{u}\|_{L^\frac{3}{2}}\le C(\|\nabla\dot{u}\|_{L^\frac{3}{2}}+\|\nabla u\|_{L^2}^{2}),
\ea\enn which together with Sobolev's embedding theorem yields \eqref{tb90}.

Finally, taking  $v=\dot{u}+(u\cdot\nabla n)\times u^{\perp}$ in \eqref{paz11} proves  \eqref{tb11} due to \eqref{paz2}.
\end{proof}

\section{\la{se3} A priori estimates (I): lower order estimates}

   Let $T>0$ be a fixed time and $(\n,u)$ be a smooth solution to (\ref{a1})-(\ref{ch1})  on
$\Omega \times (0,T]$  with smooth initial
data $(\n_0,u_0)$ satisfying (\ref{dt1}) and (\ref{dt2}). In this section, we always assume that $\Omega$ is a simply connected bounded domain in $\r^3$ and its boundary $\partial\Omega$ is of class $C^{\infty}$ and only has a finite number of 2-dimensional connected components. Since $u\cdot n=0$ on $\partial\Omega$, by \eqref{paz1}, for any $1<q<+\infty$, there exists a positive constants $C_1$ depending only on $q$, $\mu$, $\lambda$ and $\Omega$ such that
\be\la{tdu1}\ba
\|\nabla u\|_{L^q}\leq C_1(\|\div u\|_{L^q}+\|\curl u\|_{L^q}).
\ea\ee

We will establish some necessary a priori bounds for smooth solutions to the problem (\ref{a1})-(\ref{ch1}) to extend the local  classical solutions guaranteed by
Lemma \ref{loc1}.

 For $\si=\si(t)\triangleq\min\{1,t \} $  and $\dot u$ as in \eqref{oio1}, we define
 \be\la{As1}
  A_1(T) \triangleq \sup_{   0\le t\le T  }\left(\sigma\|\nabla u\|_{L^2}^2\right) + \int_0^{T} \int\sigma
 \n|\dot{u} |^2 dxdt,
  \ee
\be \la{As2}
  A_2(T)  \triangleq\sup_{  0\le t\le T   }\sigma^3\int\n|\dot{u}|^2dx + \int_0^{T}\int
  \sigma^3|\nabla\dot{u}|^2dxdt,
\ee
and
\be \la{As3}
  A_3(T)  \triangleq\sup_{  0\le t\le T   }\int\rho|u|^3dx .
\ee

Now we will give the main result in this section, which guarantees the existence of a global classical solution of \eqref{a1}--\eqref{ch1}. 
\begin{proposition}\la{pr1}  Under  the conditions of Theorem \ref{th1}, for $\delta_0\triangleq\frac{2s-1}{4s}\in(0,\frac{1}{4}],$
   there exists a  positive constant  $\ve$
    depending    on  $\mu$, $\lambda$, $a$, $\ga$,  $\hat{\rho}$, $s$, $\Omega$, $M$ and the matrix $A$ such that if
       $(\rho,u)$  is a smooth solution of
       \eqref{a1}--\eqref{ch1}  on $\Omega\times (0,T] $
        satisfying
 \be\la{zz1}
 \sup\limits_{
 \Omega\times [0,T]}\rho\le 2\hat{\rho},\quad
     A_1(T) + A_2(T) \le 2C_0^{1/3},\quad A_3(\sigma(T))\leq 2C_0^{\delta_0},
  \ee
 then the following estimates hold
        \be\la{zz2}
 \sup\limits_{\Omega\times [0,T]}\rho\le 7\hat{\rho}/4, \quad
     A_1(T) + A_2(T) \le  C_0^{1/3},\quad A_3(\sigma(T))\leq C_0^{\delta_0},
  \ee
   provided $C_0\le \ve.$
\end{proposition}
\begin{proof}Proposition \ref{pr1} is a consequence of the
following Lemmas \ref{nzc1}--\ref{le7}.
\end{proof}


One can extend the function $n$ to $\Omega$ such that $n\in C^3(\bar{\Omega})$, and in the following discussion we still denote the extended function by $n$.

In the following, we will use the convention that $C$ denotes a generic positive constant depending on $\mu,\,\,\lambda,\,\,\ga,\,\,a,\,\,\hat{\rho},\,\, s,\,\,\Omega$, $M$ and the matrix $A$, and use $C(\alpha)$ to emphasize that $C$ depends on $\alpha$.

We begin with the following standard energy estimate for $(\rho,u)$.




\begin{lemma}\la{le2}
 Let $(\n,u)$ be a smooth solution of
 \eqref{a1}--\eqref{ch1} on $\O \times (0,T]  .$
  Then there is a positive constant
  $C $ depending only  on $\mu$, $\lambda$ and $\Omega$  such that for $G(\n)$ as in \eqref{pn01},
\be \la{a16} \sup_{0\le t\le T} \int
\left( \n |u|^2+G(\n)\right)dx + \int_0^{T} \|\na u\|_{L^2}^2  dt\le C C_0 .\ee
\end{lemma}

\begin{proof}
First,
 integrating $\eqref{a1}_1$ over $\O\times (0,T)$ and using \eqref{ch1}, one has \be \la{bz11}\bar\n=\frac{1}{|\O|}\int\n (x,t)dx\equiv \frac{1}{|\O|} \int \n_0dx. \ee

 Next, since $$-\Delta u=-\nabla\div u+\nabla\times\curl u,$$ we
rewrite $(\ref{a1})_2 $ as
\be\la{m1} \ba
\rho \dot{u} - (\lambda + 2\mu)\nabla\div u+\mu\nabla\times\curl u + \nabla P  =0.
\ea \ee
Multiplying \eqref{m1} by $u$ and integrating the resulting equality over $\Omega$, along with \eqref{ch1}, we arrive at
\be\la{m8} \ba
&\frac12\left(\int\rho |u|^{2}dx\right)_t + (\lambda + 2\mu)\int(\div u)^{2}dx + \mu\int|\curl u|^{2}dx\\& +\mu\int_{\partial\Omega} u\cdot A\cdot uds=   \int P \div udx.
\ea \ee

Finally, by $(\ref{a1})_1$, one can check that
 \bnn  \ba
 (G(\n))_t+\div(G(\n)u)+(P-P(\bar\n))\div u=0,
 \ea \enn
which together with \eqref{m8} gives
\be \la{eeq1} \ba
\left(\int \frac{1}{2}\rho |u|^2+G(\rho)dx\right)_t+\phi(t)=0,
\ea\ee with
$$\phi(t)\triangleq (\lambda+2\mu)\|\div u \|_{L^{2}}^{2}+\mu\|\curl u\|_{L^{2}}^{2}+\mu\int_{\partial\Omega} u\cdot A\cdot uds.$$
  The energy estimate \eqref{a16} thus follows from  the positive semi-definiteness of $A$, \eqref{m8} and \eqref{tdu1}.
\end{proof}

A direct consequence of Lemma \ref{le2} is the following estimates on the pressure.
\begin{lemma}   \la{th01}
 Let $(\n,u)$ be a smooth solution of
 \eqref{a1}-\eqref{ch1} on $\Omega\times(0,T] $ satisfying \eqref{zz1}.
  Then there is a positive constant
  $C $ depending only  on $\mu$, $\lambda$, $a$, $\gamma$, $\hat{\rho}$, $\Omega$ and the matrix $A$ such that
\be\la{e5}\sigma\int(P-\Bar P)^{2}dx  \leq CC_0^{\frac{1}{2}}, \quad \int_0^T\si\int(P-\Bar P)^{2}dxdt\leq CC_0^{\frac{3}{4}}.\ee
 \end{lemma}
\begin{proof}  First, it follows from
 $(\ref{a1})_1$ that
 \be\la{Pu1} \ba
 P_t+\div(Pu)+(\gamma-1)P\div u=0,
 \ea \ee
which gives \be \la{ny10}\bar P_t+(\ga-1)\overline{P\div u}=0.\ee

Next, multiplying $\eqref{a1}_2$ by $\mathcal{B}[P-\Bar P]$ and integrating the resulting equality over $\Omega $, we get
\be\la{e4} \ba
\int(P-\Bar P)^2 dx &= \left(\int\rho u\cdot\mathcal{B}[P-\Bar P] dx\right)_t-\int\rho u\cdot\nabla\mathcal{B}[P-\Bar P]\cdot udx  \\
& \quad - \int\rho u\cdot\mathcal{B}[P_t-\Bar P_t]  dx +\mu\int\nabla u\cdot\nabla\mathcal{B}[P-\Bar P] dx \\
& \quad  + (\lambda+\mu)\int(P-\Bar P)\div udx \\
& \leq \left(\int\rho u\cdot\mathcal{B}[P-\Bar P] dx\right)_t+C\| u\|_{L^6}^{2}\|P-\Bar P\|_{L^{3/2}}\\
& \quad  +C\|  u\|_{L^2}\|\nabla u\|_{L^2} +C\|P-\Bar P\|_{L^2}\|\nabla u\|_{L^2}\\
& \leq \left(\int\rho u\cdot\mathcal{B}[P-\Bar P] dx\right)_t+\de \|P-\Bar P\|_{L^2}^2 +C(\de)\|\na u\|_{L^2}^2,
\ea\ee where in the first inequality we have used \bnn\ba \|\mathcal{B}[P_t-\bar P_t]\|_{L^2}&=\|\mathcal{B} [\div(Pu)] + (\ga-1) \mathcal{B} [P\div u-\ol{P\div u}]\|_{L^2} \\&\le C\|\na u\|_{L^2}.\ea\enn Combining \eqref{e4}, \eqref{a16}, and Lemma \ref{th00}  shows
\be\la{e6}\int_0^T\int(P-\Bar P)^{2}dxdt\leq CC_0^{\frac{1}{2}}.\ee

Next, using \eqref{a1}$_1,$ we have
 \be\la{Pu2} \ba
 P_t+u\cdot\nabla P +\gamma P\div u=0,
 \ea \ee  which together with   \eqref{ny10} gives
\be\la{e7}(P-\Bar P)_t+ u\cdot\na(P-\Bar P) + \gamma  P \div u -(\gamma-1)\overline{  P \div u}=0.\ee
Multiplying \eqref{e7} by $2\sigma(P-\Bar P)$ and integrating the result over $\Omega$, one checks that
\bnn \left(\sigma \int (P-\Bar P)^2 dx\right)_t
 \leq C (\si+\si')\int (P-\Bar P)^2 dx+C\int |\nabla u|^2 dx ,\enn
which together with \eqref{a16} and \eqref{e6} leads to \be  \la{ui100}  \sigma \int (P-\Bar P)^2 dx \le CC_0^{1/2}.\ee

Finally, combining \eqref{ui100} with \eqref{a16}  implies  \bnn \ba\sigma \left|\int\rho u\cdot\mathcal{B}[P-\Bar P] dx\right| &\le C\left(\int\n|u|^2dx\right)^{1/2} \left(\sigma \int (P-\Bar P)^2 dx\right)^{1/2}\\&\le CC_0^{3/4},\ea\enn which together with \eqref{a16}, \eqref{e4} and \eqref{e6} yields
\bnn\int_0^T
\si\int(P-\Bar P)^{2}dxdt\leq CC_0^{\frac{3}{4 }}.\enn Combining this with \eqref{ui100}   proves Lemma \ref{th01}.
\end{proof}

The following conclusion concerns   preliminary estimates on the $L^{2}$-norm of $\nabla u$ and $\rho^{1/2}\dot{u}$.
\begin{lemma}\la{xcrle1}
 Let $(\n,u)$ be a smooth solution of
 \eqref{a1}-\eqref{ch1} on $\Omega\times(0,T] $ satisfying \eqref{zz1}.
  Then there is a positive constant
  $C $ depending only  on $\mu$, $\lambda$, $a$, $\gamma$, $\hat{\rho}$, $\Omega$ and the matrix $A$ such that
  \be\la{h14}
  A_1(T) \le  C C_0^\frac{1}{2} + C\int_0^{T}\sigma\int|\nabla u|^3dx dt,
  \ee
 and
  \be\la{h15}
    A_2(T)
    \le   C C_0 + CA_1(T)  + C\int_0^{T}\sigma^3\int  |\nabla u|^4 dxdt.
   \ee
\end{lemma}

\begin{proof}
  Motivated by Hoff \cite{H3},
multiplying $(\ref{a1})_2 $ by
$\sigma^m \dot{u}$ with $m\ge 0 $   and then integrating the resulting equality over
$\Omega$, one gets
\be\la{I0} \ba  \int \sigma^m \rho|\dot{u}|^2dx &
= -\int\sigma^m \dot{u}\cdot\nabla Pdx + (\lambda+2\mu)\int\sigma^m \nabla\div u\cdot\dot{u}dx \\
&\quad - \mu\int\sigma^m \nabla\times\curl u\cdot\dot{u}dx \\
& \triangleq I_1+I_2+I_3. \ea \ee
We will estimate $I_1$, $I_2$ and $I_3$ one by one. First, it follows from \eqref{Pu1} that
\be\la{I10} \ba
I_1 = & - \int \sigma^m \dot{u}\cdot\nabla P dx  \\
= & \int\sigma^m P \,\div u_{t}\,dx - \int\sigma^mu\cdot\nabla u\cdot\nabla P dx  \\
= & \left(\int\sigma^m P   \div u  dx\right)_{t} - m\sigma^{m-1}\sigma'\int  P \div u\,dx + \int\sigma^{m}P\nabla u:\nabla u dx \\
&+ (\gamma-1)\int\sigma^{m}P(\div u)^{2}dx  - \int_{\partial\Omega}\sigma^{m}Pu\cdot\nabla u\cdot n ds\\
\le & \left(\int\sigma^m P \div u  dx\right)_{t} +C\|\na u\|_{L^2}^2+Cm\sigma^{m-1}\sigma' \|P-\bar P\|_{L^2}^2 ,\ea \ee
where in the last inequality we have used
\bnn\la{bdt1} \ba
- \int_{\partial\Omega}\sigma^{m}Pu\cdot\nabla u\cdot n ds&=\int_{\partial\Omega}\sigma^{m}Pu\cdot\nabla n\cdot uds
\\&\leq C\int_{\partial\Omega}\sigma^{m}|u|^{2}ds
  \leq C\sigma^{m}\|\nabla u\|_{L^{2}}^{2},
\ea  \enn due to  \eqref{pzw1}.
Hence,
\be\la{I1} \ba
I_1 \leq &\left(\int\sigma^m P\,\div u\, dx\right)_{t} + C \|\nabla u\|_{L^{2}}^{2} + Cm\sigma^{m-1}\sigma' \|P-\bar P\|^2_{L^2}.
\ea \ee

Similarly, by \eqref{pzw1}, it indicates that
\be \la{bz13}\ba
\frac{I_2}{\lambda+2\mu} & =  \int\sigma^m \nabla\div u\cdot\dot{u}dx \\
& =  \int_{\partial\Omega}\sigma^m\div u\,(\dot{u}\cdot n)ds - \int\sigma^m\div u\,\div \dot{u}dx  \\
& =  \int_{\partial\Omega}\sigma^m\div u\,(u\cdot\nabla u\cdot n)ds - \frac{1}{2}\left(\int\sigma^{m}(\div u)^{2}dx\right)_{t} \\
&\quad -  \int\sigma^m\div u\,\div(u\cdot\nabla u)dx + \frac{m }{2}\sigma^{m-1}\sigma'\int(\div u)^{2}dx \\
&\le  - \int_{\partial\Omega}\sigma^m\div u\,(u\cdot\nabla n\cdot u)ds - \frac{1}{2}\left(\int\sigma^{m}(\div u)^{2}dx\right)_{t} \\
&\quad +C\int\sigma^{m}|
\na u|^{3}dx  + \frac{m }{2}\sigma^{m-1}\sigma'\int(\div u)^{2}dx.
\ea  \ee
For the first term on the righthand side of  \eqref{bz13},  we have
\bnn \ba
&\left|(\lambda +2\mu)\int_{\partial\Omega}\div u\,(u\cdot\nabla n\cdot u)ds\right| \\
&= \left|\int_{\partial\Omega}(F+P-\Bar P)(u\cdot\nabla n\cdot u)ds\right|  \\
& \leq C(\| F\|_{H^1} + 1) \|u\|_{H^1}^{2} \\
& \leq\frac{1}{4}\|\rho^{\frac{1}{2}}\dot{u}\|_{L^{2}}^{2}+C(\|\nabla u\|_{L^{2}}^{2}+\|\nabla u\|_{L^{2}}^{4}),
\ea  \enn where in the last inequality, we have used \be \la{bz5} \|F\|_{H^1}+\|\curl u\|_{H^1}\le C(\|\n \dot u\|_{L^2}+\|\na  u\|_{L^2})\ee due to   Lemma \ref{le3}  and \eqref{tb90}. Therefore,
\be\la{I2} \ba
I_2 & \leq - \frac{\lambda+2\mu}{2}\left(\int\sigma^{m}(\div u)^{2}dx\right)_t +C\sigma^{m}\|\nabla u\|_{L^{3}}^{3}\\
&\quad +\frac{1}{4}\sigma^{m}\|\rho^{\frac{1}{2}}\dot{u}\|_{L^{2}}^{2}+C\sigma^{m}\|\nabla u\|_{L^{2}}^{4}+C\|\nabla u\|_{L^{2}}^{2}.\ea\ee

Finally, by \eqref{ch1} and \eqref{h19}, a straightforward computation shows that
\be\la{I3}\ba
I_3 & = -\mu\int\sigma^{m}\nabla\times\curl u\cdot\dot{u}dx \\
& =  - \mu\int\sigma^{m}\curl u\cdot\curl\dot{u}dx +\mu\int_{\partial\Omega}\sigma^{m}\curl u\times n\cdot\dot{u}ds \\
& = -\frac{\mu}{2}\left(\int\sigma^{m}|\curl u|^{2}dx +\int_{\partial\Omega}\sigma^{m}u\cdot A\cdot uds \right)_t \\
& \quad + \frac{\mu m}{2}\sigma^{m-1}\sigma'\int|\curl u|^{2}dx  + \frac{\mu m}{2}\sigma^{m-1}\sigma'\int_{\partial\Omega}u\cdot A\cdot uds  \\
& \quad- \mu\int\sigma^{m}\curl u\cdot\curl(u\cdot\nabla u)dx  +\mu\int_{\partial\Omega}\sigma^{m}(Au)^\perp\times(u\cdot\nabla u)\cdot nds\\
& \leq -\frac{\mu}{2}\left(\int\sigma^{m}|\curl u|^{2}dx +\int_{\partial\Omega}\sigma^{m}u\cdot A\cdot uds \right)_t + C \|\nabla u\|_{L^{2}}^{2}\\
&\quad + C\sigma^{m}\|\nabla u\|_{L^{3}}^{3}  +C \sigma^{m}\|\nabla u\|_{L^{2}}^{4}+\frac{1}{4}\sigma^{m}\|\rho^{\frac{1}{2}}\dot{u}\|_{L^{2}}^{2},
\ea \ee where in the last inequality, we have used
 \bnn\ba &\mu\int_{\partial\Omega} (Au)^\perp\times(u\cdot\nabla u)\cdot nds\\
& =     \mu\int \curl((Au)^\perp)\cdot(u\cdot\nabla u)dx - \mu\int (\nabla u^{i}\times\nabla_i u)\cdot(Au)^\perp dx \\
&\quad -\mu\int (u\cdot\nabla\curl u )\cdot(Au)^\perp dx\\&\le   C \|\nabla u\|_{L^{3}}^{3}  +C(\eta) \|\nabla u\|_{L^{2}}^{4}+\eta \|\rho^{\frac{1}{2}}\dot{u}\|_{L^{2}}^{2},\ea\enn due to \eqref{bz5}.
 It thus follows from  \eqref{I0}   and  \eqref{I1}-\eqref{I3} that
\be\la{I4}\ba
&\left(\sigma^{m}\int\left( (2\mu+\lm)(\div u)^{2}+\mu |\curl u|^{2}\right)dx +\mu\sigma^{m}\int_{\partial\Omega}u\cdot A\cdot uds \right)_t\\&\quad+\int\sigma^{m}\rho|\dot{u}|^{2}dx \\
& \leq \left(2\int\sigma^{m}(P-\bar P)\,\div udx\right)_{t}+Cm\sigma^{m-1}\sigma' \|P-\bar P\|^2_{L^2}\\
& \quad+C\sigma^{m}\|\nabla u\|_{L^{2}}^{4} +C\|\nabla u\|_{L^{2}}^{2}+C\sigma^{m}\|\nabla u\|_{L^{3}}^{3} .
\ea \ee
Integrating \eqref{I4} over $(0,T)$, by \eqref{tdu1}, \eqref{a16},   \eqref{e4},  \eqref{e6} and Young's inequality, we conclude that for any $m\ge 1$,
\be\la{I40}\ba
&\sigma^{m}\|\nabla u\|_{L^{2}}^{2}+\int_0^T\int\sigma^{m}\rho|\dot{u}|^{2}dxdt \\
& \leq CC_{0}^\frac{1}{2}+C\int_0^T\sigma^{m}\|\nabla u\|_{L^{2}}^{4}dt+C\int_0^T\sigma^{m}\|\nabla u\|_{L^{3}}^{3}dt,
\ea \ee which, after choosing $m=1,$  together with
   \eqref{zz1} and \eqref{a16} gives $\eqref{h14}$.

Now we will prove \eqref{h15}. First,  rewrite $ (\ref{a1})_2 $ as
\be\la{xdy1}\ba
\rho\dot{u}=\nabla F - \mu\nabla\times\curl u.
\ea \ee
Operating $ \sigma^{m}\dot u^j[\pa/\pa t+\div
(u\cdot)] $ to $ (\ref{xdy1})^j,$ summing with respect to $j$, and integrating it over $\Omega$, one has
\be\la{ax1}\ba &\left(\frac{\sigma^{m}}{2}\int\rho|\dot{u}|^{2}dx\right)_t -\frac{m}{2}\sigma^{m-1}\sigma'\int\rho|\dot{u}|^{2}dx \\
& = \int\sigma^{m}(\dot{u}\cdot\nabla F_t+\dot{u}^{j}\,\div(u\partial_jF))dx\\
&\quad+\mu\int\sigma^{m}(-\dot{u}\cdot\nabla\times\curl u_t-\dot{u}^{j}\div((\nabla\times\curl u)^j\,u))dx \\
& \triangleq J_1+\mu J_2.
\ea\ee

For $J_1,$ it follows from  \eqref{ny10} that
\bnn\ba F_t&=(2\mu+\lm)\div u_t-P_t+\bar P_t \\&=(2\mu+\lm)\div\dot  u-(2\mu+\lm)\na u:\na u  - u\cdot\na F+\ga P\div u-(\gamma-1)\overline{P\div u},\ea\enn which together with  \eqref{ch1}, \eqref{Pu2},  \eqref{bz5} and \eqref{bz6} yields
\be\la{ax2}\ba J_1 & =\int\sigma^{m}\dot{u}\cdot\nabla F_tdx+\int\sigma^{m}\dot{u}^{j}\div(u\partial_jF)dx \\
& = \int_{\partial\Omega}\sigma^{m}F_t\dot{u}\cdot nds-\int\sigma^{m}F_t\,\div\dot{u}dx- \int\sigma^{m}u\cdot \na \dot u^j\partial_jF   dx \\
& = \int_{\partial\Omega}\sigma^{m}F_t\dot{u}\cdot nds-(2\mu+\lm)\int\sigma^{m}(\div\dot  u)^2dx\\&\quad+(2\mu+\lm)\int\sigma^{m} \div\dot{u}\na u:\na u dx +\int\sigma^{m} \div\dot{u} u\cdot\na Fdx\\&\quad-\ga \int\sigma^{m} P\div u \div\dot{u}dx- \int\sigma^{m}u\cdot \na \dot u^j\partial_jF   dx\\&\quad
-(\gamma-1)\overline{P\div u}\int\sigma^{m}\div \dot{u} dx\\
& \le  \int_{\partial\Omega}\sigma^{m}F_t\dot{u}\cdot nds-(2\mu+\lm)\int\sigma^{m}(\div\dot  u)^2dx\\
&\quad+C\left(\|\na u\|_{L^4}^2+\|\na u\|_{L^2}+\|\na u\|_{L^2}\|\na F\|_{L^2}^\frac{1}{2}\|\na F\|_{L^6}^\frac{1}{2}\right)\|\na\dot u\|_{L^2}\\
& \le  \int_{\partial\Omega}\sigma^{m}F_t\dot{u}\cdot nds-(2\mu+\lm)\int\sigma^{m}(\div\dot  u)^2dx+\de\si^m \|\na\dot u\|_{L^2}^2\\&\quad+C(\de)\si^m\left(\|\na u\|_{L^2}^4 \|\rho^{\frac{1}{2}}\dot{u}\|_{L^2}^2+\|\na u\|_{L^4}^4+\|\na u\|_{L^2}^2\right),\ea\ee
 where in the last  equality we have used \eqref{bz5} and
\be \la{bz6}\ba & \|\na F\|_{L^6}+\|\na\curl u\|_{L^6}\\  &\le C\| \dot u\|_{H^1}+ C\|\na  u\|_{L^2}+C\|P-\bar P\|_{L^6} \\&\le C(  \|\na \dot u\|_{L^2} +\|\na u\|_{L^2}+\|\na u\|^2_{L^2}+\|P-\bar P\|_{L^6}),\ea \ee  due to Lemma \ref{le3}  and \eqref{tb90}.

For the first  term on the righthand side of \eqref{ax2}, we have
\be\la{bz8}\ba
&\int_{\partial\Omega}\sigma^{m}F_t\dot{u}\cdot nds\\&=-\int_{\partial\Omega}\sigma^{m}F_t\,(u\cdot\nabla n\cdot u)ds \\
& = -\left(\int_{\partial\Omega}\sigma^{m}(u\cdot\nabla n\cdot u)Fds\right)_t+m\sigma^{m-1}\sigma'\int_{\partial\Omega}(u\cdot\nabla n\cdot u)Fds\\
&\quad+\sigma^{m} \int_{\partial\Omega}F \dot{u}\cdot\na n\cdot uds +\sigma^{m} \int_{\partial\Omega}F {u}\cdot\na n\cdot \dot uds \\
&\quad  -\sigma^{m} \int_{\partial\Omega}F  ({u}\cdot\na)u\cdot\na n\cdot uds  -\sigma^{m} \int_{\partial\Omega}F  u\cdot\na n\cdot ({u}\cdot\na) uds  \\
& \le  -\left(\int_{\partial\Omega}\sigma^{m}(u\cdot\nabla n\cdot u)Fds\right)_t+Cm\si'\si^{m-1}\|\na u\|_{L^2}^2\|F\|_{H^1}\\&\quad +\de\si^m \|\dot u\|_{H^1}^2+C(\de)\si^{m }\|\na u\|_{L^2}^2\|F\|^2_{H^1}  \\
&\quad -\sigma^{m} \int_{\partial\Omega}F  ({u}\cdot\na)u\cdot\na n\cdot uds-\sigma^{m} \int_{\partial\Omega}F  u\cdot\na n\cdot ({u}\cdot\na) uds, \ea\ee
 where in the last inequality we have used \be \la{Fnn1} \left|\int_{\partial\Omega}(u\cdot\nabla n\cdot u)Fds \right|\le C\|\na u\|_{L^2}^2\|F\|_{H^1}.\ee

Now we are in a position to estimate the last term on the righthand side of \eqref{bz8} which indeed plays a crucial role in our analysis.  Since $u\cdot n|_{\pa\O}=0,$  we observe that \bnn  \la{bz2}  u=-(u\times n)\times n\triangleq  u^\bot\times n \mbox{ on } \pa\O,\enn which yields that \be \la{bz3}\ba &- \int_{\partial\Omega}F  ({u}\cdot\na)u\cdot\na n\cdot uds \\&= -\int_{\partial\Omega}   u^\bot\times n \cdot\na u^i \nabla_i n\cdot u Fds \\&= - \int_{\partial\Omega}n\cdot ( \na u^i \times  u^\bot)    \nabla_i n\cdot u Fds\\&= - \int_{ \Omega}\div( ( \na u^i \times  u^\bot)   \nabla_i n\cdot uF)dx \\&= - \int_{ \Omega}\na (\nabla_i n\cdot uF) \cdot ( \na u^i \times  u^\bot)   dx  + \int_{ \Omega}   \na u^i \cdot \na\times  u^\bot     \nabla_i n\cdot u  F     dx \\& \le C \int_\O |\na F||\na u||u|^2dx+C \int_\O |F| (|\na u|^2|u|+|\na u||u|^2)dx
\\& \le C  \|\na F\|_{L^6}\|\na u\|_{L^2}\|u\|^2_{L^6}+C  \| F\|_{L^{12/5}}\|\na u\|_{L^4}\|u\|^2_{L^6}  \\&\quad +C  \| F\|_{L^3}\|\na u\|^2_{L^4}\|u\|_{L^6}
\\& \le C  \|\na F\|_{L^6} \|\na u\|^3_{L^2}+ C  \|  F\|_{H^1}\|\na u\|_{L^2} \left(\|\na u\|^2_{L^4} +\|\na u\|^2_{L^2}\right) \\& \le \de \|\nabla\dot{u}\|_{L^2}^2+C(\de) \|\na u\|^6_{L^2}+C (\|\na u\|^2_{L^2}+\|\na u\|^4_{L^4})\\&\quad+ C  \|\rho^{\frac{1}{2}}\dot{u}\|_{L^2}^2 \left( \|\na u\|^2_{L^2}+1\right), \ea\ee
where in the fourth equality we have used the following standard fact:$$\div(   \na u^i \times  u^\bot)=-\na u^i\cdot \na\times u^\bot. $$

Similarly, we get
\be \la{bz4}\ba &-  \int_{\partial\Omega}F  u\cdot\na n\cdot ({u}\cdot\na) uds\\& \le \de \|\nabla\dot{u}\|_{L^2}^2+C(\de) \|\na u\|^6_{L^2}+C (\|\na u\|^2_{L^2}+\|\na u\|^4_{L^4})\\&\quad+ C  \|\rho^{\frac{1}{2}}\dot{u}\|_{L^2}^2 \left( \|\na u\|^2_{L^2}+1\right).  \ea\ee
Substituting \eqref{bz8}, \eqref{bz3}, and \eqref{bz4} into \eqref{ax2}, we obtain  after   using    \eqref{tb90}, \eqref{a16}, \eqref{bz5}, and \eqref{bz6} that
\be\la{ax399}\ba
 J_1 &\leq  Cm\sigma^{m-1}\sigma'(\|\rho^{\frac{1}{2}}\dot{u}\|_{L^2}^2+\|\nabla u\|_{L^2}^2+\|\nabla u\|_{L^2}^4)\\
&\quad-\left(\int_{\partial\Omega}\sigma^{m}(u\cdot\nabla n\cdot u)Fds\right)_t +C\delta\sigma^{m}\|\nabla\dot{u}\|_{L^2}^2\\
&\quad +C(\delta)\sigma^{m}\|\rho^{\frac{1}{2}}\dot{u}\|_{L^2}^2(\|\nabla u\|_{L^2}^4+1)- (\lambda+2\mu)\int\sigma^{m}(\div\dot{u})^{2}dx\\
& \quad+C(\delta)\sigma^{m}(\|\nabla u\|_{L^2}^2+\|\nabla u\|_{L^2}^6+\|\nabla u\|_{L^4}^4).
\ea\ee

For $J_2,$ since $\curl u_t=\curl \dot u-u\cdot \na \curl u-\na u^i\times \nabla_iu$,
\be\la{ax3999}\ba
J_2
&=- \int\sigma^{m}|\curl\dot{u}|^{2}dx+\int\sigma^{m}\curl\dot{u}\cdot(\nabla u^i\times\nabla_i u) dx\\&\quad +\int_{\partial\Omega}\sigma^{m}\curl u_t\times n\cdot\dot{u}ds +\int\sigma^{m} u\cdot\na  \curl u \cdot\curl\dot{u} dx  \\&\quad+\int\sigma^{m}    u \cdot \na \dot{u}\cdot (\nabla\times\curl u) dx  \\
&\le - \int\sigma^{m}|\curl\dot{u}|^{2}dx+\de  \sigma^{m}\|\na \dot u\|_{L^2}^2 +C(\de)(\|\na u\|^2_{L^2}+\|\na u\|^2_{L^6})\\&\quad +C(\de) \sigma^{m} \|\na u\|_{L^4}^4 +C(\de)  \sigma^{m}\|\na  u\|_{L^2}^4\|\na \curl u\|_{L^2}^2,\ea\ee
where in the last inequality we have utilized \eqref{bz5}, \eqref{bz6} and the fact
\be\la{ax4000}\ba
 \int_{\partial\Omega}\curl u_t\times n\cdot\dot{u}ds&=-\int_{\partial\Omega}u_t\cdot A\cdot\dot{u}ds\\
&=-\int_{\partial\Omega}\dot{u}\cdot A\cdot\dot{u}ds+\int_{\partial\Omega}(u\cdot\nabla u)\cdot A\cdot\dot{u}ds\\
&=-\int_{\partial\Omega}\dot{u}\cdot A\cdot\dot{u}ds+\int_{\partial\Omega}   u^\bot\times n \cdot\na u^i\, (A_{i.}\cdot \dot{u}) ds\\
&=-\int_{\partial\Omega}\dot{u}\cdot A\cdot\dot{u}ds+
\int_{\partial\Omega}n\cdot ( \na u^i \times  u^\bot)    A_{i.}\cdot \dot{u}ds\\
&=-\int_{\partial\Omega}\dot{u}\cdot A\cdot\dot{u}ds +\int_{ \Omega}\div( ( \na u^i \times  u^\bot)   A_{i.} \cdot \dot{u})dx \\
&=  -\int_{\partial\Omega}\dot{u}\cdot A\cdot\dot{u}ds+\int_{ \Omega}\na (A_{i.} \cdot \dot{u}) \cdot ( \na u^i \times  u^\bot)   dx \\&\quad - \int_{ \Omega}   \na u^i \cdot \na\times  u^\bot     (A_{i.} \cdot \dot{u})   dx,
\ea\ee
here the symbol $A_{i.}$ denotes the $i$-th row of the matrix $A$.

Putting \eqref{ax399} and \eqref{ax3999} into \eqref{ax1} yields
\be\la{ax40}\ba
&\left(\frac{\sigma^{m}}{2}\|\rho^{\frac{1}{2}}\dot{u}\|_{L^2}^2\right)_t+(\lambda+2\mu)\sigma^{m}\|\div\dot{u}\|_{L^2}^2+\mu\sigma^{m}\|\curl\dot{u}\|_{L^2}^2+\int_{\partial\Omega}\dot{u}\cdot A\cdot\dot{u}ds\\
& \leq Cm\sigma^{m-1}\sigma'(\|\rho^{\frac{1}{2}}\dot{u}\|_{L^2}^2+\|\nabla u\|_{L^2}^2+\|\nabla u\|_{L^2}^4)-\left(\int_{\partial\Omega}\sigma^{m}(u\cdot\nabla n\cdot u)Fds\right)_t \\
&\quad+\delta\sigma^{m}\|\nabla\dot{u}\|_{L^2}^2+C\sigma^{m}\|\rho^{\frac{1}{2}}\dot{u}\|_{L^2}^2(\|\nabla u\|_{L^2}^4+1)\\
& \quad+C(\delta)\sigma^{m}(\|\nabla u\|_{L^2}^2+\|\nabla u\|_{L^2}^6+\|\nabla u\|_{L^4}^4),
\ea\ee which together with
 \eqref{tb11},  after choosing $\delta$ suitably small, implies
\be\la{ax401}\ba
&\left(\sigma^{m}\|\rho^{\frac{1}{2}}\dot{u}\|_{L^2}^2\right)_t +\mu\Lambda^{-1}\sigma^{m}\|\na\dot{u}\|_{L^2}^2\\
& \leq Cm\sigma^{m-1}\sigma'(\|\rho^{\frac{1}{2}}\dot{u}\|_{L^2}^2+\|\nabla u\|_{L^2}^2+\|\nabla u\|_{L^2}^4) \\
&\quad-\left(2\int_{\partial\Omega}\sigma^{m}(u\cdot\nabla n\cdot u)Fds\right)_t+C\sigma^{m}\|\rho^{\frac{1}{2}}\dot{u}\|_{L^2}^2(\|\nabla u\|_{L^2}^4+1) \\
&\quad+C\sigma^{m}(\|\nabla u\|_{L^2}^2+\|\nabla u\|_{L^2}^6+\|\nabla u\|_{L^4}^4).
\ea\ee

Finally, integrating  \eqref{ax401} with  $m=3$ over $(0,T)$, we get \eqref{h15} from  \eqref{Fnn1} and \eqref{zz1}, which completes the proof of Lemma \ref{xcrle1}.
\end{proof}

\begin{lemma}\la{zc1} Assume that $(\n,u)$ is a smooth solution of
 \eqref{a1}-\eqref{ch1} satisfying \eqref{zz1}.  Then there exist  positive constants $C$ and $\varepsilon_1$ depending only on $\mu,\,\,\lambda,\,\,\gamma,\,\,a,\,\,\hat{\rho}, \,\, s,\,\,\Omega$, $M$ and the matrix $A$ such that
   \be\la{uv1}  \sup_{0\le t\le \si(T)}t^{1-s}\|\na
u\|_{L^2}^2+\int_0^{\si(T)}t^{1-s}\int\n|\dot u|^2dxdt\le
C(\hat{\rho},M), \ee
 \be\la{uv2}  \sup_{0\le t\le \si(T)}t^{2-s}\int\n|\dot u|^2dx+\int_0^{\si(T)}t^{2-s}\int|\nabla\dot{u}|^2dxdt\le
C(\hat{\rho},M), \ee
provide that $C_0\leq \varepsilon_1$.
\end{lemma}

\begin{proof} First, for $Lf\triangleq\rho\dot{f}-\mu\Delta f-(\lambda+\mu)\nabla\div f ,$ suppose that $w_1(x,t)$ and $w_2(x,t)$ solve the following problems  respectively
\be\la{fcz1}\begin{cases}
  Lw_1=0, \,\, &x\in\Omega,\\  w_1(x,0)=w_{10}(x),\,\, &x\in\Omega,\\w_1\cdot n=0,\,\curl w_1\times n=-Aw_1 \,\,&x\in\partial\Omega,
\end{cases}\ee
and
\be\la{fcz2}\begin{cases}
  Lw_2=-\nabla(P-\Bar P),\,\, &x\in\Omega,\\  w_2(x,0)=0,\,\, &x\in\Omega,\\w_2\cdot n=0,\,\curl w_2\times n=-Aw_2, \,\,&x\in\partial\Omega.
\end{cases}\ee

Then, just as what we have done in the proof of Lemma \ref{le3}, by Lemma \ref{zhle} and Sobolev's inequality, for any $p\in[2,6],$ we have
$$\|\nabla^{2} w_1\|_{L^{2}}\leq C(\|\rho\dot{w}_1\|_{L^{2}}+\|\nabla w_1\|_{L^{2}}),$$
\be \la{xbh1} \ba
\|\nabla w_1\|_{L^{p}}\leq C\|w_1\|_{W^{2,2}}\leq C(\|\rho\dot{w}_1\|_{L^{2}}+\|\nabla w_1\|_{L^{2}}),
\ea\ee
\be \la{xbh2} \ba
\|\nabla F_2\|_{L^{p}}\leq C(\|\rho\dot{w}_2\|_{L^{p}}+\|\nabla w_2\|_{L^{2}}+\|P-\Bar P\|_{L^{p}}),
\ea\ee
\be \la{xbh3} \ba
\|F_2\|_{L^{p}}&\leq C\|\nabla F_2\|_{L^{2}}\le C (\|\rho\dot{w}_2\|_{L^{2}}+\|\nabla w_2\|_{L^{2}}),
\ea\ee
\be \la{xbh4} \ba
\|\nabla w_2\|_{L^{p}}&\leq C\|\rho^{\frac{1}{2}}\dot{w}_2\|_{L^{2}}^{\frac{3p-6}{2p}}(\|\nabla w_2\|_{L^{2}}+\|P-\Bar P\|_{L^{2}})^{\frac{6-p}{2p}}\\
& \quad +C(\|\nabla w_2\|_{L^{2}}+\|P-\Bar P\|_{L^{p}}),
\ea\ee
where we denote the viscous effective flux of $w_2$ as $$F_2\triangleq(\lambda+2\mu)\div w_2-(P-\Bar P)  .$$

Next, a similar way as for the proof of \eqref{a16} shows that
\be\la{xbh5}  \sup_{0\le t\le \si(T)}\int\rho|w_1|^2dx+\int_0^{\si(T)}\int|\nabla w_1|^2dxdt\le
C\int|w_{10}|^{2}dx  ,\ee
and
\be\la{xbh6}  \sup_{0\le t\le \si(T)}\int\rho|w_2|^2dx+\int_0^{\si(T)}\int|\nabla w_2|^2dxdt\le
CC_0  .\ee

Then, multiplying \eqref{fcz1} by $w_{1t}$ and integrating over $\Omega,$ by \eqref{xbh1}, \eqref{zz1}, and Young's inequality, we obtain
\be \la{xbh9} \ba
&\left(\frac{\lambda+2\mu}{2}\int(\div w_1)^{2}dx + \frac{\mu}{2}\int|\curl w_1|^{2}dx +\frac{\mu}{2}\int_{\partial\Omega}w_1\cdot A\cdot w_1ds \right)_t + \int\rho|\dot{w}_1|^{2}dx \\
& =\int\rho\dot{w}_1\cdot(u\cdot\nabla w_1)dx \\
& \le C\|\rho^{\frac{1}{2}}\dot{w}_1\|_{L^{2}}\|\rho^{\frac{1}{3}}u\|_{L^{3}}\|\nabla w_1\|_{L^{6}} \\
& \le CC_0^{\frac{\delta_0}{3}}(\|\rho^{\frac{1}{2}}\dot{w}_1\|_{L^{2}}^{2}+ \|\nabla w_1\|_{L^{2}}^{2}) ,
\ea\ee
which together with \eqref{xbh5},  Gronwall's inequality and Lemma \ref{crle1}  yields
\be \la{xbh10} \ba
\sup_{0\le t\le \si(T)}\|\nabla w_1\|_{L^{2}}^{2}+\int_0^{\sigma(T)}\int\rho|\dot{w}_1|^{2}dxdt\leq C\|\nabla w_{10}\|_{L^{2}}^{2}  ,
\ea\ee
and
\be \la{xbh11} \ba
\sup_{0\le t\le \si(T)}t\|\nabla w_1\|_{L^{2}}^{2}+\int_0^{\sigma(T)}t\int\rho|\dot{w}_1|^{2}dxdt\leq C\|w_{10}\|_{L^{2}}^{2},
\ea\ee
provided $C_0\leq \hat{\varepsilon}_1\triangleq(2C)^{-\frac{3}{\delta_0}}$.

Next, since the solution operator $w_{10}\mapsto w_1(\cdot,t)$ is linear,  one can deduce from Calder\'{o}n's interpolation theorem \cite[Lemmas 22.3 and 36.1]{TL1},  (\ref{xbh10}) and (\ref{xbh11}) that for any $\theta\in [s,1],$
\be \la{xbh12} \ba
 \sup_{0\le t\le
\si(T)}t^{1-\theta}\|\nabla
w_1\|_{L^2}^2+\int_0^{\si(T)}t^{1-\theta}\int\n|\dot
{w_1}|^2dxdt\leq C\| w_{10}\|_{H^\theta}^2,
\ea\ee
with a uniform constant $C$ independent of $\theta.$

Next, multiplying (\ref{fcz2}) by $w_{2t}$ and integrating over $\Omega$, we give
\be \la{xbh13} \ba
&\left(\frac{\lambda+2\mu}{2}\int(\div w_2)^{2}dx+\frac{\mu}{2}\int|\curl w_2|^{2}dx-\int P\div w_2dx  \right)_t\\
&\quad
+\left( \frac{\mu}{2}\int_{\partial\Omega}w_2\cdot A\cdot w_2ds \right)_t+\int\rho|\dot{w}_2|^{2}dx \\
& =\int\rho\dot{w}_2\cdot(u\cdot\nabla w_2)dx-\int P_t\div w_2dx\\
& =\int\rho\dot{w}_2\cdot(u\cdot\nabla w_2)dx - \frac{1}{\lambda+2\mu}\int P(F_2\div u+\nabla F_2\cdot u ) dx\\
& \quad - \frac{1}{2(\lambda+2\mu)}\int(P-\Bar P)^{2}\div u dx +\gamma\int P\div u\,\div w_2dx \\
& \leq C(\|\rho^{\frac{1}{2}}\dot{w}_2\|_{L^2}\|\rho^{\frac{1}{3}}u\|_{L^3}\|\nabla w_2\|_{L^6}+\|\nabla u\|_{L^2}\|F_2\|_{L^2}+\|\nabla F_2\|_{L^2}\|u\|_{L^2})\\
&\quad+C(\|P-\Bar P\|_{L^2}\|\nabla u\|_{L^2}+\|\nabla u\|_{L^2}\|\nabla w_2\|_{L^2})\\
&\leq CC_0^{\frac{\delta_0}{3}}\|\rho^{\frac{1}{2}}\dot{w}_2\|_{L^2}(\|\rho^{\frac{1}{2}}\dot{w}_2\|_{L^2}+\|\nabla w_2\|_{L^2}+\|P-\Bar P\|_{L^6})\\
&\quad+C\|\nabla u\|_{L^2}\|\rho^{\frac{1}{2}}\dot{w}_2\|_{L^2}+C(\|P-\Bar P\|_{L^2}\|\nabla u\|_{L^2}+\|\nabla u\|_{L^2}\|\nabla w_2\|_{L^2})\\
&\leq CC_0^{\frac{\delta_0}{3}}\|\rho^{\frac{1}{2}}\dot{w}_2\|_{L^2}^{2}+\frac{1}{4}\|\rho^{\frac{1}{2}}\dot{w}_2\|_{L^2}^{2}
+C(\|\nabla w_2\|_{L^2}^2+\|\nabla u\|_{L^2}^2+\|P-\Bar P\|_{L^2}^\frac{2}{3}),
\ea\ee
where we have utilized \eqref{zz1}, \eqref{Pu2}, \eqref{xbh2}-\eqref{xbh4}, H\"{o}lder's, Poincar\'{e}'s and Young's inequalities.  After choosing  $C_0\leq \hat{\varepsilon}_2\triangleq(4C)^{-\frac{3}{\delta_0}} ,$  combining this  with   Gronwall's inequality, \eqref{xbh6}, and Lemmas \ref{crle1}, \ref{le2} and \ref{th01}, we get
\be \la{xbh15} \ba
\sup_{0\le t\le \si(T)}\|\nabla w_2\|_{L^{2}}^{2}+\int_0^{\sigma(T)}\int\rho|\dot{w}_2|^{2}dxdt\leq C .
\ea\ee

Now letting $w_{10}=u_0 $, we have $w_1+w_2=u,$ which combined with  (\ref{xbh12}) and (\ref{xbh15}) directly proves
  \eqref{uv1}  provided $C_0\leq \varepsilon_1\triangleq\min\{\hat{\varepsilon}_1,\hat{\varepsilon}_2\}$.

Finally, it remains   to prove \eqref{uv2}.  Taking $m=2-s$ in \eqref{ax401},  and integrating over $(0,\sigma(T)),$  we obtain by \eqref{tb11},
\be \la{xbh16} \ba
&\sup_{0\le t\le  \si(T)}t^{2-s}\|\rho^{\frac{1}{2}}\dot{u}\|_{L^2}^2+\int_0^{\sigma(T)}t^{2-s}\|\nabla\dot{u}\|_{L^2}^2dt\\
& \leq C\int_0^{\sigma(T)}t^{1-s}\|\rho^{\frac{1}{2}}\dot{u}\|_{L^2}^2dt
+C\int_0^{\sigma(T)}t^{2-s}\|\rho^{\frac{1}{2}}\dot{u}\|_{L^2}^2(\|\nabla u\|_{L^2}^4+1)dt\\
&\quad +C\int_0^{\sigma(T)}t^{2-s}(\|\nabla u\|_{L^2}^2+\|\nabla u\|_{L^2}^6)dt+C\int_0^{\sigma(T)}t^{2-s}\|\nabla u\|_{L^4}^4dt\\
& \quad +C\int_0^{\sigma(T)}t^{1-s}(\|\nabla u\|_{L^2}^2+\|\nabla u\|_{L^2}^4)dt+Ct^{2-s}(\|\nabla u\|_{L^2}^2+\|\nabla u\|_{L^2}^4)\\
& \le  C\int_0^{\si(T)}t^{2-s}\|\na u\|_{L^4}^4dt+C(\hat{\rho}, M), \ea\ee
where we have taken advantage of \eqref{Fnn1} and \eqref{uv1}.

By \eqref{h18} and \eqref{uv1}, we have
\bnn \la{xbh17}\ba
&\int_0^{\sigma(T)}t^{2-s}\|\nabla u\|_{L^{4}}^{4}dt \\
& \le C\int_0^{\sigma(T)}t^{2-s}\|\rho^{\frac{1}{2}}\dot{u}\|_{L^{2}}^{3}(\|\nabla u\|_{L^{2}}+\|P-\Bar P\|_{L^{2}})dt\\
& \quad +C\int_0^{\sigma(T)}t^{2-s}(\|\nabla u\|_{L^{2}}^{4}+\|P-\Bar P\|_{L^{4}}^{4})dt \\
& \le C\int_0^{\sigma(T)}t^{\frac{2s-1}{2}}(t^{1-s}\|\nabla u\|_{L^2}^2)^{\frac{1}{2}}(t^{2-s}\|\rho^{\frac{1}{2}}\dot{u}\|_{L^{2}}^{2})^{\frac{1}{2}}(t^{1-s}\|\rho^{\frac{1}{2}}\dot{u}\|_{L^2}^2)dt+C\\
& \leq C(\hat{\rho}, M)\left(\sup_{0\le t\le  \si(T)}t^{2-s}\|\rho^{\frac{1}{2}}\dot{u}\|_{L^2}^2\right)^{\frac{1}{2}}+C ,
\ea\enn
which together  with  \eqref{xbh16} gives \eqref{uv2}.
\end{proof}
\begin{lemma}\la{nzc1} If $(\n,u)$ is a smooth solution of \eqref{a1}-\eqref{ch1} satisfying \eqref{zz1} and the initial data condition $\|u_0\|_{H^s}\leq M$ in \eqref{dt2}, then there
exists a positive constant  $\varepsilon_2$   depending only on $\mu,\,\,\lambda,\,\,\ga,\,\,a,\,\,\hat{\rho},\,\,s,\,\,\Omega$, $M$ and the matrix $A$ such
that
\be\la{xuv1} \sup_{0\le t\le  \si(T) }\int \n |u|^{3}dx\le C_0^{\delta_0} ,\ee
provided $C_0\leq \varepsilon_2$.
\end{lemma}

\begin{proof}
First, multiplying $(\ref{a1})_2$ by $3|u|u$  and integrating the resulting equality over $ \O$, we find that
\bnn\ba & \left(\int \n |u|^{3}dx\right)_t+3(\lambda+2\mu)\int\div u\,\div(|u|u)dx+ 3\mu\int\curl u\cdot\curl(|u|u)dx  \\
& \quad  +3\mu\int_{\partial\Omega}|u|u\cdot A\cdot uds  - 3\int(P-\Bar P)\div(|u|u)dx =0,\ea\enn
which together with   \eqref{uv1} and  \eqref{uv2} yields
\be\la{ki01}\ba
&\left(\int\rho|u|^{3}dx\right)_t\\
&\leq C\int|u||\nabla u|^{2}dx+C\int_{\pa\O}|u|^3ds+C\int|P-\Bar P||u||\nabla u|dx\\
&\leq C\|u\|_{L^6}\|\nabla u\|_{L^2}^{\frac{3}{2}}\|\nabla u\|_{L^6}^{\frac{1}{2}}+C\|\na u\|_{L^2}^3+C\|P-\Bar P\|_{L^3}\|u\|_{L^6}\|\nabla u\|_{L^2}\\
&\leq C\|\nabla u\|_{L^2}^{\frac{5}{2}}\si^{-\frac{2-s}{4}}+C\|\nabla u\|_{L^2}^3 +C\|\nabla u\|_{L^2}^2\\
&\leq C\|\nabla u\|_{L^2}^{4\de_0}\si^{-\frac{(6-8\de_0)(1-s)+1}{4}} ,
\ea\ee
where in the third inequality we have used
\bnn\ba \|\na u\|_{L^6}&\le C\|\rho\dot{u}\|_{L^2}+C\|P-\Bar P\|_{L^2}+C\|P-\Bar P\|_{L^6}+C\|\nabla u\|_{L^2}\\&\le C\si^{-(2-s)/2}\ea\enn
due to \eqref{h18}.

Thus, combining \eqref{ki01} and   \eqref{a16} implies
\bnn\la{bcbh1} \ba
&\sup_{0\le t\le  \si(T) }\int\rho|u|^{3}dx\\
&\leq C(\hat{\rho}, M)\left(\int_0^{\sigma(T)}\si^{-\frac{(6-8\de_0)(1-s)+1}{4(1-2\de_0)}}dt \right)^{ 1-2\delta_0}  \left(\int_0^{\sigma(T)}\|\nabla u\|_{L^2}^2dt\right)^{2\delta_0 }\\
&\quad +\int\rho_0|u_0|^3dx\\
&\leq C(\hat{\rho}, M)C_0^{2\delta_0},
\ea\enn
provided $C_0\leq \varepsilon_1,$ where in the last inequality we have used
both $\frac{(6-8\de_0)(1-s)+1}{4(1-2\de_0)}<1 $
due to $\delta_0=\frac{2s-1}{4s}\in(0,\frac{1}{4}]$ and $s\in(\frac{1}{2},1]$ and the following simple fact (see\cite[Theorem 1]{BM1})
\be\la{bcbh2} \ba
& \int\rho_0|u_0|^{3}dx\leq C\|\rho_0^{\frac{1}{2}}u_0\|_{L^{2}}^{{3(2s-1)}/{2s}}\|u_0\|_{H^s}^{{3}/{2s}}\leq C(\hat{\rho}, M)C_0^{2\delta_0}.
\ea\ee

Finally, we obtain \eqref{xuv1} by setting $\varepsilon_2\triangleq\min\{\varepsilon_1,(C(\hat{\rho}, M))^{-\frac{1}{\delta_0}}\}.$ The proof of Lemma \ref{nzc1} is finished.
\end{proof}

\begin{lemma}\la{le5} Let $(\n,u)$ be a smooth solution  of
   \eqref{a1}-\eqref{ch1}     on $\O \times (0,T] $ satisfying \eqref{zz1} and the initial data condition $\|u_0\|_{H^s}\leq M$ in \eqref{dt2}. Then there exists a positive constant $\varepsilon_3$ depending only  on $\mu$,  $\lambda$,  $\gamma$, $a$, $s$, $\hat{\rho}$, $M,$ $\Omega$ and the matrix $A$
 such that
  \be\la{h27}
  A_1(T)+A_2(T)\le C_0^{\frac{1}{3}},
  \ee
 provided $C_0\leq\varepsilon_3$.
   \end{lemma}

\begin{proof}  First, we will prove \eqref{h27}. By (\ref{h18}) and \eqref{zz1},  one can check that
  \be\la{h9h9} \ba
   \sigma^3 \|\na u\|_{L^4}^4
& \le  C  \sigma^{3} \|\n^{\frac{1}{2}}  \dot u \|_{L^2}^3(\|\nabla u\|_{L^{2}}+\|P-\Bar P\|_{L^{2}}) \\
& \quad + C \sigma^{3}(\|\nabla u\|_{L^{2}}^{4}+\|P-\Bar P\|_{L^{4}}^{4}) \\
& \le C C^{1/3}_0\sigma  \|\n^{\frac{1}{2}}  \dot u \|_{L^2}^2+C(\|\nabla u\|^2_{L^{2}}+\|P-\Bar P\|^2_{L^{2}}),
 \ea \ee which together with \eqref{zz1}, \eqref{a16} and \eqref{e5} leads to
  \be\la{h99} \ba
  &  \int_0^{T}\sigma^3 \|\na u\|_{L^4}^4 dt\le  CC_0^\frac{1}{2}.
 \ea \ee
Next, it follows from (\ref{h18}), \eqref{zz1}, \eqref{a16}, \eqref{uv1} and \eqref{e5} that
\be\la{h992} \ba
&\int_0^{\sigma(T)}\sigma\|\nabla u\|_{L^{3}}^{3}dt \\
& \le C\int_0^{\sigma(T)}\sigma\|\rho^{\frac{1}{2}}\dot{u}\|_{L^{2}}^{\frac{3}{2}}(\|\nabla u\|_{L^{2}}^{\frac{3}{2}}+\|P-\Bar P\|_{L^{2}}^{\frac{3}{2}})dt\\
& \quad +C\int_0^{\sigma(T)}\sigma(\|\nabla u\|_{L^{2}}^{3}+\|P-\Bar P\|_{L^{3}}^{3})dt\\
& \le C\int_0^{\sigma(T)}(\si^{\frac{1-s}{2}}\|\nabla u\|_{L^{2}})\|\nabla u\|_{L^{2}}^{\frac{1}{2}}(\sigma\|\rho^{\frac{1}{2}}\dot{u}\|_{L^{2}}^{2})^{\frac{3}{4}}dt\\
& \quad +C\int_0^{\sigma(T)}\sigma^{\frac{1}{4}}\|P-\Bar P\|_{L^{2}}^{\frac{3}{2}}(\sigma\|\rho^{\frac{1}{2}}\dot{u}\|_{L^{2}}^{2})^{\frac{3}{4}}dt\\
& \quad +C\left(\int_0^{\sigma(T)}(\sigma\|\nabla u\|_{L^{2}}^{2})\|\nabla u\|_{L^{2}}dt+\int_0^{\sigma(T)}\sigma\|P-\Bar P\|_{L^{2}}^{2}dt\right)\\
\quad &\le C(M)\left(\int_0^{\sigma(T)}\|\nabla u\|_{L^{2}}^{2}dt\right)^{\frac{1}{4}}\left(\int_0^{\sigma(T)}\sigma\|\rho\dot{u}\|_{L^{2}}^{2}dt\right)^{\frac{3}{4}}+CC_0^{\frac{1}{2}} \\
& \le C(\hat{\rho}, M)(A_1(T))^{\frac{3}{4}}C_0^{\frac{1}{8}}+CC_0^{\frac{1}{2}}\\
& \le C(\hat{\rho}, M)C_0^{\frac{3}{8}},
\ea \ee
provided $C_0\leq \varepsilon_2$.

On the other hand, by \eqref{h99} and \eqref{a16},
\be\la{h993} \ba
 \int_{\sigma(T)}^T\sigma\|\nabla u\|_{L^{3}}^{3}dt
& \le\int_{\sigma(T)}^T\sigma\|\nabla u\|_{L^{4}}^{4}dt+\int_{\sigma(T)}^T\sigma\|\nabla u\|_{L^{2}}^{2}dt\\
& \le CC_0^\frac{1}{2} .
\ea \ee

Finally, by \eqref{h14}, \eqref{h15} and \eqref{h99}-\eqref{h993}, we have
\bnn\la{h994} \ba
A_1(T)+A_2(T)\leq C(\hat{\rho}, M)C_0^{\frac{3}{8}},
\ea \enn which
 gives \eqref{h27} provided  $C_0\leq \varepsilon_3$ with $\varepsilon_3\triangleq\min\{\varepsilon_2,(C(\hat{\rho}, M))^{-24}\}.$
\end{proof}

We now proceed to derive a uniform (in time) upper bound for the
density, which turns out to be the key to obtaining all the higher
order estimates and thus to extending the classical solution globally in time.
We will use an approach motivated by the works  \cite{lx,hlx1}.

\begin{lemma}\la{le7}
There exists a positive constant
   $\ve$
    depending    on  $\mu$,  $\lambda$, $\ga$, $a$, $\hat{\rho}$, $s$, $ \Omega,$ $M$, and the matrix $A$  such that,
    if  $(\n,u)$ is a smooth solution  of
   \eqref{a1}-\eqref{ch1}     on $\O \times (0,T] $
   satisfying \eqref{zz1} and the initial data condition $\|u_0\|_{H^s}\leq M$ in \eqref{dt2}, then for  $(x,t)\in \O \times (0,T)$
      \be\la{lv102}\sup_{0\le t\le T}\|\n(t)\|_{L^\infty}  \le
\frac{7\hat \n }{4}  ,\ee
      provided $C_0\le \ve. $ Moreover, if  $C_0\le \ve , $ there exists some positive constant $\ti C(T)$ depending only on $T,$ $\mu$,  $\lambda$, $\ga$, $a$, $\hat{\rho}$, $s$, $ \Omega,$ $M$, and the matrix $A$  such that  for  $(x,t)\in \O \times (0,T)$
      \be \la{lvv79} \n(x,t)\ge \ti C(T)\inf_{x\in\O}\n_0(x). \ee
\end{lemma}

\begin{proof}
  First, the equation of  mass conservation $(\ref{a1})_1$ can be  rewritten in the form
 \be \la{z.3} D_t \n=g(\n)+b'(t), \ee where \bnn
D_t\n\triangleq\n_t+u \cdot\nabla \n ,\quad
g(\n)\triangleq-\frac{\rho (P- \Bar P)}{2\mu+\lambda},
\quad b(t)\triangleq -\frac{1}{2\mu+\lambda} \int_0^t  \rho F dt. \enn


Then, it follows from \eqref{bz6},  \eqref{zz1},  \eqref{uv1}, and \eqref{uv2} that for $\delta_0$ as in Proposition \ref{pr1}  and $t\in [0,\si(T)],$
 \bnn \|\na F\|_{L^6}\le C\|\na \dot u\|_{L^2}+C\si^{-(1-s)},\,\,\|\rho^{\frac{1}{2}}\dot{u}\|_{L^2} \le C \si^{-\frac{(2-s)(1-\de_0)+3\de_0}{2}}C_0^{\frac{ \de_0}{6}},\enn
provide $C_0\le \varepsilon_1$. Combining this,  \eqref{g2},  \eqref{hh22},   and \eqref{zz1}   yields
\bnn  \ba
 \|  F (\cdot,t)\|_{L^{\infty}}
& \le C \|F\|_{L^{6}}^{\frac{1}{2}}\|\nabla F\|_{L^{6}}^{\frac{1}{2}}   \\
& \le C \left( \|\rho^{\frac{1}{2}}\dot{u}\|_{L^2}  +\|\nabla u\|_{L^2}\right)^{\frac{1}{2}}\left(\|\na\dot{u}\|_{L^{2}}+\si^{-(1-s)} \right)^{\frac{1}{2}}    \\
& \le C   \si^{-\frac{(2-s)(2-\de_0)+3\de_0 }{4}}C_0^{\frac{ \de_0}{ 12}}\left(\si^{2-s}\|\na\dot{u}\|^2_{L^{2}}+\si^{s} \right)^{\frac{1}{4}},
\ea\enn
 which together with  \eqref{uv2}  and Holder's inequality thus implies that
 for all $0\leq t_1\leq t_2\leq\sigma(T)$,
\be \la{xbh19} \ba
 |b(t_2)-b(t_1)|  &
\le C\int_0^{\sigma(T)}\|  F (\cdot,t)\|_{L^{\infty}}dt
  \le C C_0^{\frac{ \de_0}{ 12}}  ,
\ea\ee  due to $(2-s)(2-\de_0)+3\de_0<3.$
 Thus, choosing $N_1=0$, $N_0=C(\hat{\rho}, M)C_0^{\frac{\delta_0}{12}}$, and $\zeta_0=\hat{\rho}$
in Lemma  \ref{le1}, we use   \eqref{xbh19}, \eqref{z.3}, and Lemma  \ref{le1} to get
   \be\la{a103}\sup_{t\in
[0,\si(T)]}\|\n\|_{L^\infty} \le \hat{\rho}
+C(\hat{\rho}, M)C_0^{\frac{\delta_0}{12}}\le\frac{3 \hat{\n}  }{2},\ee
 provided $C_0\le \ve_4\triangleq\min\left\{\varepsilon_3, \left(\frac{\hat{\rho}}{2C(\hat{\rho}, M)}\right)^{\frac{12}{\delta_0}}\right\}. $

On the other hand, for $\sigma(T)\le t_1\le t_2\le T ,$  we have
\be \la{xbh20} \ba
 |b(t_2)-b(t_1)|
&\le C\int_{t_1}^{t_2}\|F\|_{L^{\infty}}dt\\
&\le \frac{a\hat{\rho}^{\gamma+1}}{2(\lambda+2\mu)}(t_2-t_1)+C\int_{\sigma(T)}^{T} \|F\|_{L^{\infty}}^4dt  \\&\le \frac{a\hat{\rho}^{\gamma+1}}{2(\lambda+2\mu)}(t_2-t_1)+CC_0^{\frac{1}{2}},
\ea\ee
where in the last inequality we have used
\be \la{xs20} \ba
 \int_{\sigma(T)}^{T} \|F\|_{L^{\infty}}^4dt
&\le C\int_{\sigma(T)}^{T} \|F\|_{L^6}^2\|\nabla F\|_{L^6}^{2}dt \\
& \le   CC_0^{\frac{1}{3}}\int_{\sigma(T)}^{T}\|\nabla \dot{u}\|_{L^2}^{2}dt+CC_0^{\frac{1}{2}} \\&\le  CC_0^{\frac{1}{2}},
\ea\ee due to  \eqref{h19}, \eqref{hh22}, \eqref{tb90}, \eqref{zz1}, \eqref{a16} and \eqref{e5}.

Now  choosing $N_0=CC_0^{\frac{1}{2}}$, $N_1=\frac{a\hat{\rho}^{\gamma+1}}{2(\lambda+2\mu)}$ in \eqref{a100} and setting $\zeta_0=\frac{3\hat{\rho}}{2}$ in (\ref{a101}),  we have  for all $  \zeta \geq\zeta_0=\frac{3\hat{\rho}}{2}$,
$$ g(\zeta)=-\frac{ \zeta(a\zeta^{\gamma}-\bar P)}{\lambda+2\mu}\le -\frac{a\hat{\rho}^{\gamma+1}}{2(\lambda+2\mu)}= -N_1, $$ which together with
    Lemma \ref{le1},   \eqref{a103}, and   \eqref{xbh20} leads to
\be\la{a102} \sup_{t\in
[\si(T),T]}\|\n\|_{L^\infty}\le  \frac{3\hat \n }{2} +CC_0^{\frac{1}{2}} \le
\frac{7\hat \n }{4} ,\ee provided
\be \la{xbh21} \ba C_0\le
\ve \triangleq\min\left\{ {\ve}_4, (\frac{\hat \n }{4C})^2\right\}.
\ea\ee
 Combining \eqref{a101} and \eqref{a102} thus gives \eqref{lv102}.

Finally, it remains to prove \eqref{lvv79}. Indeed, without loss of generality, assume that $\inf\limits_{x\in \O}\n_0(x)>0.$  We have by \eqref{z.3},
\bnn\ba (2\mu+\lm)D_t\n^{-1}-\n^{-1}(P-\bar P)-\n^{-1}F=0,\ea\enn
which in particular shows that
\bnn\ba  D_t\n^{-1} &\le C \n^{-1}\left(\|F\|_{L^\infty}+1\right) . \ea\enn
Combining this with Gronwall's inequality, \eqref{xbh19} and \eqref{xs20}  gives \eqref{lvv79} and finishes the proof of Lemma \ref{le7}.
\end{proof}

With Proposition \ref{pr1} at  hand, we are now in a position to prove
the following result concerning the exponential decay rate   of both weak and classical solutions. It should be noted here that  both the rate $\eta_0$ and the constant $C$    depend  on $\bar \n_0$ also which is different from the constants $C$ in the proof of Proposition \ref{pr1} where they are independent of $\bar\rho_0.$
\begin{proposition}\la{beha1}   For any $r\in [1,\infty)$ and $p\in [1,6],$ there exist positive constants $C$ and $\eta_0$ depending only  on $\mu,$  $\lambda,$  $\gamma,$ $a$, $s$,  $\hat{\rho}$, $\bar\n_0$, $M$,  $\Omega$, $r$, $p,$  and the matrix $A$   such that \eqref{qa1w} holds for $ t\geq1.$
\end{proposition}
 \begin{proof}
First, by \eqref{pn01} and \eqref{bz11}, there exists a positive constant $\tilde{C} <1$ depending only on $\gamma$, $\bar \rho_0 ,$ and $\hat \n$ such that for any $\rho\in [0,2\hat\n]$,
\be\label{gine1}  \tilde{C}^2( \rho-\bar{\rho})^2\le \tilde{C}  G(\rho)  \leq    (\rho^\gamma-\bar{\rho}^\gamma)( \rho - \bar{\rho}), \ee
and
\be \la{uyq01} \|P-\bar P\|_{L^2}^2\le C\|P-P(\bar\n)\|_{L^2}^2\le C\int G(\n)dx.\ee

Then, multiplying $\eqref{a1}_2$ by $\mathcal{B}[\n-\bar\n]$, one has
\bnn\la{c59} \ba &
\int(P-P(\bar\n))(\n-\bar\n) dx \\&= \left(\int\rho u\cdot\mathcal{B}[\n-\bar\n] dx\right)_t-\int\rho u\cdot\nabla\mathcal{B}[\n-\bar\n]\cdot udx - \int\rho u\cdot\mathcal{B}[\n_t]  dx \\
& \quad  +\mu\int\nabla u\cdot\nabla\mathcal{B}[\n-\bar\n] dx + (\lambda+\mu)\int(\rho-\bar{\rho})\div udx \\
& \leq \left(\int\rho u\cdot\mathcal{B}[\n-\bar\n] dx\right)_t+C\|\rho^{\frac{1}{2}}u\|_{L^{4}}^{2}\|\n-\bar\n\|_{L^2} +C\|\rho u\|_{L^2}^2\\
& \quad  +C\|\rho-\bar{\rho}\|_{L^2}\|\nabla u\|_{L^2} \\
& \leq \left(\int\rho u\cdot\mathcal{B}[\n-\bar\n] dx\right)_t+\de \|\n-\bar\n\|_{L^2}^2 +C(\de)\|\na u\|_{L^2}^2,
\ea\enn
which, along with \eqref{gine1} and \eqref{tdu1}, leads to
\be \la{gine2} \ba
a\tilde{C}\int G(\rho)dx&\leq a\int(\rho^\gamma-\bar{\rho}^\gamma)( \rho - \bar{\rho})dx\\
&\leq 2\left(\int\rho u\cdot\mathcal{B}[\n-\bar\n] dx\right)_t+\ti C_1\phi(t).
\ea\ee
 Moreover,   it follows from \eqref{gine1} that
\bnn \la{c511} \ba
\left|\int\rho u\cdot\mathcal{B}[\n-\bar\n] dx\right|\leq \ti C_2 \left(\frac{1}{2}\|\sqrt{\rho} u\|^2_{L^2}+\int G(\rho)dx\right),
\ea\enn
which gives
\be \la{c512} \ba
\frac{1}{2}\left(\frac{1}{2}\|\sqrt{\rho} u\|^2_{L^2}+\int G(\rho)dx\right)\leq W(t)\leq 2\left(\frac{1}{2}\|\sqrt{\rho} u\|^2_{L^2}+\int G(\rho)dx\right),\ea\ee
 where $$W(t)=\int \left(\frac{1}{2}\rho |u|^2+G(\rho)\right)dx-\delta_1\int\rho u\cdot\mathcal{B}[\n-\bar\n] dx,$$ with $\delta_1=\min\{(2\ti C_1)^{-1},(2\ti C_2)^{-1}\}.$

Adding \eqref{gine2} multiplied by $\de_1 $ to \eqref{eeq1}
and utilizing \bnn \int\n |u|^2dx\leq C\|\na u\|_{L^2}^2\leq C_3\phi(t),\enn we obtain for $\eta_0=\min\{\frac{a\delta_1 \tilde{C}}{4},\frac{1}{4C_3}\}$,
 \bnn W'(t)+2\eta_0W(t)\leq 0.\enn Combining this with \eqref{c512} yields that for any $t>0$,
\begin{equation}\label{c513}
\int\left(\frac{1}{2}\rho|u|^2+G(\rho)\right)dx\leq 4C_0e^{-2\eta_0t},
\end{equation} which together with  \eqref{eeq1} shows
\be \la{c514} \ba
\int_0^T\phi(t)e^{\eta_0 t} dt\leq C.
\ea\ee
Choosing $m=0$ in \eqref{I4}  along with \eqref{tdu1}, \eqref{h18}, \eqref{zz1} and \eqref{uyq01}  leads to
\be\la{c515}\ba
&\left(\phi(t)  -2\int(P-P(\bar{\rho}))\,\div udx\right)_{t}\\&+\frac{1}{2}\|\sqrt{\rho}\dot{u}\|^2_{L^2}\leq C\phi(t)+C G(\rho).
\ea \ee

Noticing that
\bnn\ba
\left|\int(P-P(\bar{\rho}))\div udx\right|\leq CG(\n) +\frac{1}{4}\phi(t),
\ea\enn
we obtain after multiplying \eqref{c515} by $e^{\eta_0 t}$ and using   \eqref{c513} and \eqref{c514}  that for any $T\geq 1$,
\be\la{c517}\ba
\sup_{1\leq t\leq T}\left(e^{ \eta_0 t}\|\nabla u\|_{L^2}^2\right) +
\int_1^T e^{\eta_0 t}\|\sqrt{\rho}\dot{u}\|^2_{L^2}dt\leq C.
\ea \ee

Finally,   a similar analysis based on  \eqref{ax401} where $m=3,$ \eqref{h9h9}, \eqref{Fnn1}, \eqref{uyq01}   and \eqref{c517} indicates that for any $t\geq 1$,
\bnn\la{c519}\ba
\|\sqrt{\rho}\dot{u}\|^2_{L^2}\leq Ce^{-\eta_0 t},
\ea \enn
which together with
  \eqref{c513}, \eqref{c517},  and \eqref{h18} gives \eqref{qa1w}  and finishes the proof of Proposition \ref{beha1}.
\end{proof}
\section{\la{se5} A priori estimates (II): higher order estimates }

Let $(\n,u)$ be a smooth solution of \eqref{a1}-\eqref{ch1}. The purpose of this section is to derive some necessary higher order estimates, which make sure that one can extend the classical solution globally in time. Here we adopt the method of the article \cite{hlx1,jx01}, and follow their work with a few modifications. We sketch it here for completeness.

In this section, we always assume that the initial energy $C_0$ satisfies (\ref{xbh21}), and the positive
constant $C $ may depend on \bnn  T,\,\, \| g\|_{L^2},  \,\,\|\na u_0\|_{H^1},\,\,
    \|\n_0 \|_{W^{2,q}}  ,  \,\, \|P(\n_0) \|_{W^{2,q}} , \,\,\enn
for $q\in(3,6)$ besides $\mu$, $\lambda$, $a$, $\ga$, $\hat{\rho}$,   $s,$  $\Omega$, $
M$ and the matrix $A$, where $g\in L^2(\Omega)$ is given by \eqref{dt3}.
\begin{lemma}\la{xle1}
 There exists a positive constant $C $ such that
\be\label{cxb2}
\sup_{0\le t\le T}\|\rho^{\frac{1}{2}}\dot{u}\|_{L^2}+\int_0^T\|\nabla\dot{u}\|_{L^2}^{2}dt\leq C,\ee
\be\label{cxb3}
\sup_{0\le t\le T}(\|\nabla\rho\|_{L^6}+\|u\|_{H^2})+\int_0^T(\|\nabla u\|_{L^\infty}+\|\nabla^{2} u\|_{L^6}^{2})dt\leq C.\ee
\end{lemma}
\begin{proof} First, Taking $s=1$ in (\ref{uv1}) along with
(\ref{h27}) gives
     \be\la{cxb4}
  \sup_{t\in[0,T]}\|\nabla u\|_{L^2}^2 + \int_0^{T}\int\rho|\dot{u}|^2dxdt
  \le C .
  \ee
Choosing $m=0$ in \eqref{ax401}, by \eqref{h18} and \eqref{tb11}, we have
\be\la{cxb5} \ba
&\left(\|\rho^{\frac{1}{2}}\dot{u}\|_{L^2}^2\right)_t+
\mu\Lambda^{-1}\|\nabla\dot{u}\|_{L^2}^2\\
& \leq -\left(2\int_{\partial\Omega}(u\cdot\nabla n\cdot u)Fds\right)_t+C\|\rho^{\frac{1}{2}}\dot{u}\|_{L^2}^2(\|\nabla u\|_{L^2}^4+1) \\
&\quad+C(\|\nabla u\|_{L^2}^2+\|\nabla u\|_{L^2}^6+\|\nabla u\|_{L^4}^4)\\
& \leq -\left(2\int_{\partial\Omega}(u\cdot\nabla n\cdot u)Fds\right)_t+C\|\rho^{\frac{1}{2}}\dot{u}\|_{L^2}^2(\|\rho^{\frac{1}{2}}\dot{u}\|_{L^2}^2+\|\nabla u\|_{L^2}^4+1) \\
&\quad+C(\|\nabla u\|_{L^2}^2+\|\nabla u\|_{L^2}^6+\|P-\Bar P\|_{L^4}^4).
 \ea\ee
By Gronwall's inequality and the compatibility condition \eqref{dt3}, we deduce \eqref{cxb2} from \eqref{cxb5}, \eqref{cxb4} and \eqref{Fnn1}.

Next,  we will follow  the proof of   \cite[Lemma 5]{hlx} to  show \eqref{cxb3}.
For $2\leq p\leq 6 ,$ $|\nabla\rho|^p$ satisfies \bnnn \ba
& (|\nabla\rho|^p)_t + \text{div}(|\nabla\rho|^pu)+ (p-1)|\nabla\rho|^p\text{div}u  \\
 &+ p|\nabla\rho|^{p-2}(\nabla\rho)^{tr} \nabla u (\nabla\rho) +
p\rho|\nabla\rho|^{p-2}\nabla\rho\cdot\nabla\text{div}u = 0,\ea
\ennn
where $(\nabla\rho)^{\rm tr}$ is the transpose of $\nabla\rho$.

Thus, taking $p=6$, by \eqref{h19}, \eqref{tb90} and \eqref{cxb4},
\be\la{cxb9}\ba
(\|\nabla\rho\|_{L^6})_t&\le C(1+\norm[L^{\infty}]{\nabla u} )\norm[L^6]{\nabla\rho} +C\|\na F\|_{L^6}\\
&\le C(1+\norm[L^{\infty}]{\nabla u} )\norm[L^6]{\nabla\rho}+C\|\rho\dot{u}\|_{L^6}\\
&\le C(1+\norm[L^{\infty}]{\nabla u} )\norm[L^6]{\nabla\rho}+C(\|\nabla\dot{u}\|_{L^2}+1).
 \ea\ee

Then, it follows from the Gagliardo-Nirenberg inequality, \eqref{tb90}, \eqref{h19}, and \eqref{cxb4} that
\bnn\la{cxb11}\ba
&\|\div u\|_{L^\infty}+\|\curl u\|_{L^\infty}\\
&\le C(\|F\|_{L^\infty}+\|P-\Bar P\|_{L^\infty})+\|\curl u\|_{L^\infty} \\
&\le C(\|F\|_{L^2}+\|\nabla F\|_{L^6}+\|\curl u\|_{L^2}+\|\nabla \curl u\|_{L^6}+\|P-\Bar P\|_{L^\infty}) \\
&\le C(\|\nabla u\|_{L^2}+\|P-\Bar P\|_{L^2}+\|\rho\dot{u}\|_{L^6}+\|P-\Bar P\|_{L^\infty}) \\
&\le C(\|\nabla\dot{u}\|_{L^2}+1),
\ea\enn
which together with  Lemma \ref{le9}   and \eqref{tb90} yields
\be\la{cxb12}\ba
\|\na u\|_{L^\infty } &\le C\left(\|{\rm div}u\|_{L^\infty }+
\|\curl u\|_{L^\infty } \right)\ln(e+\|\na^2 u\|_{L^6 }) +C\|\na
u\|_{L^2} +C \\
&\le C(1+\|\nabla\dot{u}\|_{L^2})\ln(e+\|\nabla\dot u\|_{L^2 } +\|\na \rho\|_{L^6})   \\
&\le C(1+\|\nabla\dot{u}\|_{L^2}^2)+C(1+\|\nabla\dot{u}\|_{L^2})\ln(e+\|\nabla\rho\|_{L^6}) ,
\ea\ee
where in the second inequality, we have used the fact that for any $p\in[2,6]$,
\be\la{remark1}\ba
\|\nabla^{2}u\|_{L^p}\leq C(\|\rho\dot{u}\|_{L^p}+\|\nabla P\|_{L^p}+\|\nabla u\|_{L^2}+\|P-\Bar P\|_{L^p}),
\ea\ee
which can be obtained by applying Lemma \ref{zhle} to the following system
\be\la{lames}\begin{cases}
-\mu\Delta u-(\lambda+\mu)\nabla\div u=-\rho\dot{u}-\nabla(P-\Bar P), \,\, &x\in\Omega, \\
u\cdot n=0,\,\, \,\curl u\times n=-Au,\,\,&x\in\partial\Omega.
\end{cases} \ee

Next, it follows from  \eqref{cxb12} and  \eqref{cxb9} that
\bnn\la{cxb13}\ba
&(e+\|\nabla\rho\|_{L^6})_t\\
&\leq C\left(1+\|\nabla\dot{u}\|_{L^2}^2+(1+\|\nabla\dot{u}\|_{L^2}) \ln(e+\|\nabla\rho\|_{L^6})\right)(e+\|\nabla\rho\|_{L^6}),
\ea\enn
which  together with Gronwall's inequality and \eqref{cxb2} shows that
\bnn\la{cxb15}\ba
\sup_{0\leq t\leq T}\|\nabla\rho\|_{L^{6}}\leq C.
\ea\enn
 Combining  this,  \eqref{cxb12},   \eqref{remark1}, \eqref{tb90}, \eqref{cxb2}, and \eqref{cxb4} gives \eqref{cxb3} and finishes the proof of Lemma \ref{xle1}.
\end{proof}

\begin{lemma}\la{xle2}
 There exists a positive constant $C$ such that
\be\la{cxb17}\ba
\sup_{0\le t\le T}\|\rho^{\frac{1}{2}}u_t\|_{L^2}^2 + \ia\int|\nabla u_t|^2dxdt\le C,
\ea\ee
\be\la{cxb18}\ba
\sup_{0\le t\le T}(\|{\rho }\|_{H^2} +
 \|{P }\|_{H^2})\le C.
\ea\ee
\end{lemma}
\begin{proof} By Lemma \ref{xle1}, a straightforward calculation yields that
\bnn\la{cxb19}\ba
\|\rho^{\frac{1}{2}}
u_t\|_{L^2}^2 &\le \|\rho^{\frac{1}{2}}\dot u \|_{L^2}^2+\|\n^{\frac{1}{2}}u\cdot\na u\|_{L^2}^2\\
&\le C+C\|u\|_{L^4}^2\|\nabla u\|_{L^4}^2 \\
&\le C+C\|\nabla u\|_{L^2}^2\|u\|_{H^2}^2 \\
&\le C ,
\ea\enn
and
\bnn \la{cxb20}\ba
 \int_0^T\|\nabla u_t\|_{L^2}^2dt &\le\int_0^T\|\nabla \dot
u\|_{L^2}^2dt + \int_0^T\|\nabla(u\cdot\nabla u)\|_{L^2}^2dt \\
&\le C+\int_0^T\|\nabla u\|_{L^4}^4+\|u\|_{L^\infty}^2\|\nabla^{2}u\|_{L^2}^2dt  \\
&\le C+C\int_0^T(\|\nabla^{2}u\|_{L^2}^4+\|\nabla u\|_{H^1}^{2}\|\nabla^{2}u\|_{L^2}^2)dt \\
&\le C .
\ea\enn

 It remains to prove \eqref{cxb18}.
 We deduce from \eqref{lames} and Lemma \ref{zhle} that for any $p\in[2,6]$,
\be\la{remark2}\ba
\|\nabla^{3}u\|_{L^p}\leq C(\|\rho\dot{u}\|_{W^{1,p}}+\|\nabla P\|_{W^{1,p}}+\|\nabla u\|_{L^2}+\|P-\Bar P\|_{L^p}),
\ea\ee which together with
 \eqref{Pu2}, $(\ref{a1})_1$,  \eqref{remark1},  and Lemma \ref{xle1} gives
\bnn \la{cxb22}\ba
&\frac{d}{dt}\left(\|\nabla^2P\|_{L^2}^2 +\|\nabla^2\rho\|_{L^2}^2\right)\\
&\le C(1+\|\nabla u\|_{L^\infty})(\|\nabla^2P\|_{L^2}^2 +\|\nabla^2\rho
\|_{L^2}^2) + C\|\nabla\dot{u}\|_{L^2}^2 + C.
\ea\enn
Consequently, combining this,  Gronwall's inequality,  and Lemma \ref{xle1} leads to
  \bnn \sup_{0\le t\le T} {\left(\|\nabla^2P\|_{L^2}^2
+\|\nabla^2\rho\|_{L^2}^2 \right)}\le C. \enn Thus the proof of
Lemma \ref{xle2} is finished.
\end{proof}
\begin{lemma}\la{xle3}
There exists a positive constant $C $ such that
 \be\la{cxb24}
   \sup\limits_{0\le t\le T}\left(
   \|\n_t\|_{H^1}+\|P_t\|_{H^1}\right)
    + \int_0^T\left(\|\n_{tt}\|_{L^2}^2+\|P_{tt}\|_{L^2}^2\right)dt
\le C,
  \ee
\be\la{cxb25}
   \sup\limits_{0\le t\le T}\si \|\nabla u_t\|_{L^2}^2
    + \int_0^T\si\|\rho^{\frac{1}{2}}u_{tt}\|_{L^2}^2dt
\le C.
  \ee
\end{lemma}
\begin{proof} First, it follows from (\ref{Pu2}) and Lemma \ref{xle1} that
\be \la{cxb26}
\|P_t\|_{L^2}\le
C\|u\|_{L^\infty}\|\nabla P\|_{L^2}+C\|\nabla u\|_{L^2}\le C.
\ee

Next, applying $\na$ to (\ref{Pu2}), one gets
\bnn
\nabla P_t+u\cdot\nabla\nabla P+\nabla u\cdot\nabla P+\ga \nabla P {\rm div}u+\ga P  \nabla{\rm div}u=0,
\enn
which together with   Lemmas \ref{xle1} and \ref{xle2} gives
\bnn\la{cxb27} \|\nabla P_t\|_{L^2}\le C\|u\|_{L^\infty}\|\nabla^2
P\|_{L^2}+C\|\nabla u\|_{L^3}\|\nabla P\|_{L^6}+C\|\nabla^2
u\|_{L^2}\le C.\enn
Combining this with (\ref{cxb26}), one has
\bn \la{cxb28}\sup_{0\le t\le T}\|P_t\|_{H^1}\le C.
\en

 Next, it follows from \eqref{Pu2} that $P_{tt}$ satisfies
\be\la{cxb29} P_{tt} + \gamma P_t{\rm div}u +
\gamma P{\rm div}u_t + u_t\cdot\nabla P + u\cdot\nabla P_t = 0.
\ee
Multiplying \eqref{cxb29} by $P_{tt}$ and integrating over $\Omega\times[0,T],$ we deduce from \eqref{cxb28}, Lemmas \ref{xle1} and \ref{xle2} that
\bnn \ba
 \int_0^T\|P_{tt}\|_{L^2}^2dt
& = -\int_0^T\int\gamma P_{tt}P_t\div udxdt - \int_0^T\int\gamma P_{tt}P\div u_tdxdt  \\
& \quad - \int_0^T\int P_{tt}u_t\cdot\nabla Pdxdt - \int_0^T\int P_{tt}u\cdot\nabla P_tdxdt \\
& \le C\int_0^T\|P_{tt}\|_{L^2}(\|P_{t}\|_{L^3}\|\nabla u\|_{L^6}+\|\nabla u_t\|_{L^2} )dt\\
& \quad+ C\int_0^T\|P_{tt}\|_{L^2}( \|u_t\|_{L^3}\|\nabla P\|_{L^6}+ \|u\|_{L^\infty}\|\nabla P_{t}\|_{L^2}) dt\\  & \le C\int_0^T\|P_{tt}\|_{L^2}(1+\|\nabla u_t\|_{L^2})dt\\
& \le \frac{1}{2}\int_0^T\|P_{tt}\|_{L^2}^2dt+C\int_0^T\|\nabla u_t\|_{L^2}^2dt+C\\
& \le \frac{1}{2}\int_0^T\|P_{tt}\|_{L^2}^2dt+C,
\ea\enn
where we have ulitized Sobolev's inequality. Therefore, it holds
$$\int_0^T\|P_{tt}\|_{L^2}^2dt \le C .$$
One can deal with $\n_t$ and
$\n_{tt}$ similarly. Thus, (\ref{cxb24}) is proved.

It remains to prove (\ref{cxb25}).
Since $u_t\cdot n = 0$ on $\partial\Omega$, by Lemma \ref{crle1}, we have
\be\ba\la{cxb40}
\|\nabla u_t\|_{L^2}^2\leq C H(t),
\ea\ee with
$$H(t)\triangleq(\lambda+2\mu)\int(\div u_t)^{2}dx+\mu\int|\curl u_t|^{2}dx .$$
Differentiating  $(\ref{a1})_2$  with respect to $t,$ then multiplying by
$u_{tt},$  we obtain
\be\la{cxb34} \ba
&\frac{d}{dt}\left(H(t)+ \mu\int_{\partial\Omega}u_t\cdot A\cdot u_tds \right)+2\int\rho|u_{tt}|^2dx  \\
&=\frac{d}{dt}\ti I_0
  +\int\rho_{tt}|u_t|^{2}dx + 2\int(\rho_tu\cdot\nabla u)_t\cdot u_tdx \\
&\quad-2\int\rho u_t\cdot\nabla u\cdot u_{tt}dx -2\int\rho u\cdot\nabla u_t\cdot u_{tt}dx - 2\int P_{tt}\div u_tdx \\
&\triangleq\frac{d}{dt}\ti I_0 + \sum\limits_{i=1}^5\ti I_i ,
\ea \ee
where
\be \ba \la{cxb35}
\ti I_0& \triangleq-\int_{
}\rho_t |u_t|^2 dx- 2\int_{ }\rho_t u\cdot\nabla u\cdot u_tdx+
2\int_{ }P_t {\rm div}u_tdx \\
&\le \left|\int_{ } {\rm
div}(\rho u)\,|u_t|^2dx\right|+C\norm[L^3]{\rho_t}\| u\|_{L^\infty}\|\nabla u\|_{L^2}
\norm[L^6]{u_t}\\&\quad+C\| P_t\|_{L^2}\|\nabla u_t\|_{L^2}\\
&\le C \int_{}  |u||\n u_t||\nabla u_t| dx +C\|\nabla u_t\|_{L^2} \\
&\le C\|u\|_{L^6}\|\n^{1/2} u_t\|_{L^2}^{1/2}\|u_t\|_{L^6}^{1/2}\|\nabla
u_t\|_{L^2} +C\|\nabla u_t\|_{L^2}\\
&\le C\|\nabla u\|_{L^2}\|\n^{1/2} u_t\|_{L^2}^{1/2}\|\nabla u_t\|_{L^2}^{3/2}+C\|\nabla u_t\|_{L^2}\\
&\le \frac{1}{2}H(t)+C,
\ea\ee de to $(\ref{a1})_1$, \eqref{tb90}, \eqref{cxb2}, \eqref{cxb3}, \eqref{cxb4}, \eqref{cxb17},
\eqref{cxb24}, \eqref{cxb40} and Sobolev's and Poincar\'{e}'s  inequalities.

Then, standard arguments yield that
\be \la{cxb36}\ba
|\ti I_1|&=\left|\int_{ }\rho_{tt}\, |u_t|^2 dx\right|\\
&= \left|\int_{ }\div(\rho u)_t\,|u_t|^2 dx\right|\\
&= 2\left|\int_{ }(\rho_tu + \rho u_t)\cdot\nabla u_t\cdot u_tdx\right|\\
& \le  C\left(\norm[H^1]{\rho_t}\norm[H^2]{u}
  +\norm[L^2]{\rho^{{1/2}}u_t}^{\frac{1}{2}}\|\nabla u_t\|_{L^2}^{\frac{1}{2}}\right)\|\nabla u_t\|_{L^2}^2 \\
& \le C\|\nabla u_t\|_{L^2}^4+C\|\nabla u_t\|_{L^2}^2+C\\
& \le C\|\nabla u_t\|_{L^2}^2H(t)+C\|\nabla u_t\|_{L^2}^2+C,
\ea \ee
\be \la{cxb37}\ba
|\ti I_2|&=2\left|\int_{ }\left(\rho_t u\cdot\nabla u \right)_t\cdot u_{t}dx\right|\\
&= 2\left|  \int_{ }\left(\rho_{tt} u\cdot\nabla u\cdot u_t +\rho_t
u_t\cdot\nabla u\cdot u_t+\rho_t u\cdot\nabla u_t\cdot
u_t\right)dx\right|\\
&\le\norm[L^2]{\rho_{tt}}\norm[L^3]{u\cdot\nabla u}\norm[L^6]{u_t}+\norm[L^2]{\rho_t}\|u_t\|_{L^6}^2\norm[L^6]{\nabla u} \\
&\quad+\norm[L^3]{\rho_t}\norm[L^{\infty}]{u}\norm[L^2]{\nabla u_t}\norm[L^6]{u_t}\\
& \le C\norm[L^2]{\rho_{tt}}^2 + C\norm[L^2]{\nabla u_t}^2, \ea \ee
\be\ba\la{cxb38}
|\ti I_3|+|\ti I_4|&= 2\left| \int_{ }\rho u_t\cdot\nabla
u\cdot u_{tt} dx\right| +2\left| \int_{ }\rho u\cdot\nabla u_t\cdot
u_{tt} dx\right|\\& \le   C\|\n^{1/2}u_{tt}\|_{L^2}\left(
\|u_t\|_{L^6}\|\na u\|_{L^3}+\|u\|_{L^\infty}\|\na
u_t\|_{L^2}\right) \\& \le   \norm[L^2]{\rho^{{1/2}}u_{tt}}^2 +
C\norm[L^2]{\nabla u_t}^2, \ea\ee
and
\be\ba\la{cxb39}
|\ti I_5|&=2\left|\int_{ }P_{tt}{\rm div}u_tdx\right|\\&\le C
\norm[L^2]{P_{tt}}\norm[L^2]{{\rm div}u_t}\\& \le
C\norm[L^2]{P_{tt}}^2 + C\norm[L^2]{\nabla u_t}^2.
\ea\ee

Finally, putting  \eqref{cxb36}-\eqref{cxb39} into \eqref {cxb34} gives
\bnn\la{cxb41} \ba
&\frac{d}{dt}(\sigma H(t)+ \mu\sigma\int_{\partial\Omega}u_t\cdot A\cdot u_tds -\sigma \ti I_0)+\sigma\int\rho|u_{tt}|^{2}dx \\
&\le C(1+\|\nabla u_t\|_{L^2}^2)\sigma H(t)+C(1+\norm[L^2]{\nabla u_t}^2+\|\rho_{tt}\|_{L^2}^2+\|P_{tt}\|_{L^2}^2),
\ea \enn
which together with  Gronwall's inequality, \eqref{cxb17}, \eqref{cxb24} and \eqref{cxb35} leads to
\be\la{cxb42} \ba
&\sup_{0\le t\le T}(\sigma H(t)+\mu\sigma\int_{\partial\Omega}u_t\cdot A\cdot u_tds)+\int_0^T\sigma\|\rho^{\frac{1}{2}}u_{tt}\|_{L^2}^2dt\le C .
\ea \ee
As a result, by \eqref{cxb40},
\be\la{cxb43} \ba
&\sup_{0\le t\le T}\sigma \|\nabla u_t\|_{L^2}^2+\int_0^T\sigma\|\rho^{\frac{1}{2}}u_{tt}\|_{L^2}^2dt\le C ,
\ea \ee
which finishes the proof of Lemma \ref{xle3}.
\end{proof}
\begin{lemma}\la{xle4}For  $q\in(3,6),$
there exists a positive constant $C$ such that
\be\la{cxb45}\ba\sup_{t\in[0,T]} \si \|\nabla u\|_{H^2}^2  +\ia \left(\|\nabla  u\|_{H^2}^2+\|\na^2
u\|^{p_0}_{W^{1,q}}+\si\|\na u_t\|_{H^1}^2\right)dt\le C,\ea \ee
\be\la{cxb44}\ba \sup_{t\in[0,T]}\left(\|\rho \|_{W^{2,q}} +\|P \|_{W^{2,q}}\right)\le C,\ea \ee
where $p_0=\frac{9q-6}{10q-12}\in(1,\frac{7}{6}).$
\end{lemma}
 \begin{proof} First, by Lemma \ref{xle1} and Poincar\'{e}'s and  Sobolev's inequalities, one can check that
  \bnn \label{cxb47} \ba
  \|\nabla (\n \dot u) \|_{L^2}&\le
 \||\nabla \n |\, |  u_t|  \|_{L^2}+ \|\n \,\nabla   u_t  \|_{L^2}
 + \||\nabla \n|\,| u|\,|\nabla u| \|_{L^2}\\ &\quad
 + \|\n\,|\nabla  u|^2\|_{L^2}
 + \|  \n \,|u |\,| \nabla^2 u| \|_{L^2}\\&\le
 \|\nabla \n \|_{L^3} \|  u_t  \|_{L^6}+ C\| \nabla   u_t  \|_{L^2}
 + C\| \nabla \n\|_{L^3}\| u\|_{L^\infty}\|\nabla u \|_{L^6}\\
 &\quad + C\| \nabla  u\|_{L^3}\| \nabla  u\|_{L^6}
 + C\|     u \|_{L^\infty}\| \nabla^2 u  \|_{L^2}\\ &\le C+C\| \nabla   u_t  \|_{L^2},
 \ea\enn which together with \eqref{cxb18} and Lemma \ref{xle1}  yields
\bnn\la{cxb46} \ba\|\nabla^2 u\|_{H^1} &\le C (\|\n \dot u\|_{H^1}+ \| P-\Bar P\|_{H^2}+\|u\|_{L^2})\\
 &\le C+C \|\na  u_t\|_{L^2}.
 \ea\enn
Combining this ,
(\ref{cxb3}), (\ref{cxb17}) and (\ref{cxb25}) leads to
\be\la{cxb48}
\sup\limits_{0\le
t\le T}\si\|\nabla  u\|_{H^2}^2+\ia \|\nabla  u\|_{H^2}^2dt \le
 C.\ee

Next,  it follows from \eqref{a1} and \eqref{ch1} that $u_t$ satisfies
\be\begin{cases}
  \mu\Delta u_t+(\lambda+\mu)\nabla\div u_t=(\rho\dot{u})_t+\nabla P_t \,\,\, &\text{in} \,\,\Omega,\\ u_t\cdot n=0,\,\,\,\curl u_t\times n=-Au_t\,\,&\text{on} \,\,\partial\Omega ,
\end{cases}\ee  which together with  Lemmas \ref{zhle}  and  \ref{xle1}, \eqref{cxb18} and \eqref{cxb24} shows
\be\la{cxb480}\ba
\|\na^2u_t\|_{L^2}
&\le C(\|(\rho\dot{u})_t\|_{L^2}+\|P_t\|_{H^1}+\|u_t\|_{L^2}) \\
&\le C(\|\n  u_{tt}+\n_t u_t+\n_t u\cdot\nabla u + \n u_t\cdot\nabla u+\n u\cdot\nabla u_t\|_{L^2})\\
&\quad +C(\|\nabla P_t\|_{L^2}+\|u_t\|_{L^2})+C\\
&\le C\left(\|\n  u_{tt}\|_{L^2}+ \|\n_t\|_{L^3}\|u_t\|_{L^6}+\|\n_t\|_{L^3}\| u\|_{L^\infty}\|\nabla u\|_{L^6}\right)+C\\
&\quad+C\left(\| u_t\|_{L^6}\|\nabla u\|_{L^3}+ \| u\|_{L^\infty}\|\nabla u_t\|_{L^2}+\|\nabla P_t\|_{L^2}+\|u_t\|_{L^2}\right)\\
&\le C\|\n^{\frac{1}{2}}  u_{tt}\|_{L^2} +C\|\nabla  u_t\|_{L^2}+C,
\ea \ee
Combining this, \eqref{cxb480} and \eqref{cxb25} yields
\be\la{cxb50}\ba
\int_0^T\sigma\|\nabla u_t\|_{H^1}^2dt\leq C.
\ea \ee

Next, by Sobolev's inequality, \eqref{tb90}, \eqref{cxb3}, \eqref{cxb18} and \eqref{cxb25}, we get for any $q\in (3,6)$,
\be\la{cxb53}\ba     \|\na(\n\dot u)\|_{L^q}
&\le C \|\na \n\|_{L^q}(\|\nabla\dot{u}\|_{L^q}+\|\nabla\dot{u}\|_{L^2}+\|\nabla u\|_{L^2}^2)+C\|\na\dot u \|_{L^q}\\
&\le C (\|\nabla\dot{u}\|_{L^2}+\|\nabla u\|_{L^2}^2)+C(\|\na u_t \|_{L^q}+\|\na(u\cdot \na u ) \|_{L^q})\\
&\le C (\|\nabla u_t\|_{L^2}+1)+C\|\na u_t \|_{L^2}^{\frac{6-q}{2q}}\|\nabla u_t\|_{L^6}^{\frac{3(q-2)}{2q}}\\
&\quad+C(\|u \|_{L^\infty}\|\nabla^{2}u\|_{L^q}+\|\nabla u\|_{L^{\infty}}\|\nabla u\|_{L^q})\\
&\le C\sigma^{-\frac{1}{2}}+C\|\nabla u\|_{H^2}+C\sigma^{-\frac{1}{2}}(\sigma\|\nabla u_t\|_{H^1}^2)^{\frac{3(q-2)}{4q}}+C.
\ea \ee
Integrating this over $(0,T),$ by \eqref{cxb2} and \eqref{cxb50}, we have
\be\ba\la{cxb55}
\int_0^T\|\nabla(\rho\dot{u})\|_{L^q}^{p_0}dt\leq C .
\ea\ee

On the other hand, the combination of  \eqref{Pu2} with \eqref{cxb18} gives
\be\la{cxb51}\ba
(\|\na^2 P\|_{L^q})_t\le& C \|\na u\|_{L^\infty} \|\na^2 P\|_{L^q}   +C  \|\na^2 u\|_{W^{1,q}}   \\\le& C (1+\|\na u\|_{L^\infty} )\|\na^2 P\|_{L^q}+C(1+ \|\na  u_t\|_{L^2})\\&+ C\| \na(\n
\dot u )\|_{L^{q}},
\ea\ee
 where in the last inequality we have used the  following simple fact that
\be \la{cxb52}\ba \|\na^2 u\|_{W^{1,q}}
 &\le C(\|\rho\dot{u}\|_{L^q}+\|\nabla(\rho\dot{u})\|_{L^q}+\|\nabla^{2} P\|_{L^q}+\|\nabla P\|_{L^q}\\
&\quad+\|\nabla u\|_{L^2}+\|P-\Bar P\|_{L^2}+\|P-\Bar P\|_{L^q})\\
 &\le C(1 + \|\na  u_t\|_{L^2}+ \| \na(\n\dot u )\|_{L^{q}}+\|\na^2  P\|_{L^{q}}),
 \ea\ee
 due to \eqref{remark1}, \eqref{remark2}, \eqref{cxb2} and \eqref{cxb18}.

 Hence, applying Gronwall's inequality in \eqref{cxb51}, we deduce from \eqref{cxb3}, \eqref{cxb17}  and \eqref{cxb55} that
\bnn\ba\la{cxb56}
\sup_{t\in[0,T]}\|\nabla^{2}P\|_{L^q}\leq C ,
\ea\enn
which along with \eqref{cxb17}, \eqref{cxb18}, \eqref{cxb52} and \eqref{cxb55} gives
\be\ba\la{cxb57}
\sup_{t\in[0,T]}\|P \|_{W^{2,q}}+\int_0^T\|\nabla^{2}u\|_{W^{1,q}}^{p_0}dt\leq C .
\ea\ee
Similarly, one has
\bnn\sup\limits_{0\le t\le T}\|
\n \|_{W^{2,q}} \le
 C,\enn
 which together with (\ref{cxb57})  gives (\ref{cxb44}). The proof of Lemma \ref{xle4}
is finished.
\end{proof}
\begin{lemma}\la{xle5} For   $q\in (3,6),$ there exists a positive constant $C$ such that
\be \la{cxb58}
\sup_{0\le t\le T}\si\left(\|\na u_t\|_{H^1}
 +\|\na u\|_{W^{2,q}}\right)
 +\int_{0}^T\si^2\|\nabla u_{tt}\|_{L^2}^2dt\le C.
 \ee

\end{lemma}

\begin{proof} First,  differentiating $(\ref{a1})_2$ with respect to $t$ twice implies
 \be\la{cxb59}\ba
&\n u_{ttt}+\n u\cdot\na u_{tt}-(\lambda+2\mu)\nabla{\rm div}u_{tt}+\mu\nabla\times\curl u_{tt}\\
&= 2{\rm div}(\n u)u_{tt}
+{\rm div}(\n u)_{t}u_t-2(\n u)_t\cdot\na u_t-(\n_{tt} u+2\n_t u_t)
\cdot\na u\\& \quad- \n u_{tt}\cdot\na u-\na P_{tt}.
 \ea\ee

Then, multiplying (\ref{cxb59}) by $2u_{tt}$ and   integrating over $\Omega$, we get
\be \la{cxb60}\ba
&\frac{d}{dt}\int_{ }\n
|u_{tt}|^2dx+2(\lambda+2\mu)\int_{ }(\div u_{tt})^2dx+2\mu\int_{ }|\curl u_{tt}|^2dx \\
&=-8\int_{ }  \n u^i_{tt} u\cdot\na
 u^i_{tt} dx-2\int_{ }(\n u)_t\cdot \left[\na (u_t\cdot u_{tt})+2\na
u_t\cdot u_{tt}\right]dx\\&\quad -2\int_{
}(\n_{tt}u+2\n_tu_t)\cdot\na u\cdot u_{tt}dx-2\int_{ }   \n
u_{tt}\cdot\na u\cdot  u_{tt} dx\\&\quad+2\int_{ } P_{tt}{\rm
div}u_{tt}dx\triangleq\sum_{i=1}^5K_i.
\ea\ee
Let us estimate each $K_i(i=1,\cdots,5)$ as follows.
H\"{o}lder's inequality and \eqref{cxb3} give
\be \la{cxb61} \ba |K_1|&\le
C\|\n^{1/2}u_{tt}\|_{L^2}\|\na u_{tt}\|_{L^2}\| u \|_{L^\infty}\\
&\le \de \|\na u_{tt}\|_{L^2}^2+C(\de)\|\n^{1/2}u_{tt}\|^2_{L^2} .
\ea\ee
By \eqref{cxb2}, \eqref{cxb17}, \eqref{cxb24} and \eqref{cxb25}, we conclude that
\be \la{cxb62}\ba
|K_2|&\le C\left(\|\n
u_t\|_{L^3}+\|\n_t u\|_{L^3}\right)\left(\| u_{tt}\|_{L^6}\| \na
u_t\|_{L^2}+\| \na u_{tt}\|_{L^2}\| u_t\|_{L^6}\right)\\&\le
C\left(\|\n^{1/2} u_t\|^{1/2}_{L^2}\|u_t\|^{1/2}_{L^6}+\|\n_t
\|_{L^6}\| u\|_{L^6}\right)  \| \na u_{tt}\|_{L^2}  \| \na u_{t}\|_{L^2} \\ &\le \de
\|\na u_{tt}\|_{L^2}^2+C(\de)\si^{-3/2},\ea\ee

\be  \la{cxb63}\ba |K_3|&\le C\left(\|\n_{tt}\|_{L^2}
\|u\|_{L^\infty}\|\na u\|_{L^3}+\|\n_{
t}\|_{L^6}\|u_{t}\|_{L^6}\|\na u \|_{L^2}\right)\|u_{tt}\|_{L^6} \\
&\le \de \|\na u_{tt}\|_{L^2}^2+C(\de)\|\n_{tt}\|_{L^2}^2+C(\de)\si^{-1},
\ea\ee
and
\be  \la{cxb64}\ba
|K_4|+|K_5|&\le C\|\n u_{tt}\|_{L^2} \|\na
u\|_{L^3}\|u_{tt}\|_{L^6} +C \|P_{tt}\|_{L^2}\|\na
u_{tt}\|_{L^2}\\
&\le \de \|\na u_{tt}\|_{L^2}^2+C(\de)\|\n^{1/2}u_{tt}\|^2_{L^2}
+C(\de)\|P_{tt}\|^2_{L^2}.
\ea\ee
 Since $u_{tt}\cdot n=0$ on $\partial\Omega,$ Lemma \ref{crle1} implies that there exists some positive constant $\ti\mu$ depending only on $ \O $ such that
\be  \la{cxb65}\ba
\ti \mu\|\nabla u_{tt}\|_{L^2}\leq   \|\div u_{tt}\|_{L^2}+\|\curl u_{tt}\|_{L^2}  .
\ea\ee
 Thus, substituting \eqref{cxb61}--\eqref{cxb64} into (\ref{cxb60})  and  choosing $\de$ small enough, one has
\bnn  \la{cxb66}\ba
&\frac{d}{dt}\|\n^{1/2}u_{tt}\|^2_{L^2}+  \mu\ti\mu \|\na u_{tt}\|_{L^2}^2\\
&\le C (\|\n^{1/2}u_{tt}\|^2_{L^2}+\|\n_{tt}\|^2_{L^2}+\|P_{tt}\|^2_{L^2})+C \si^{-3/2},
 \ea\enn
which together with  (\ref{cxb24}), (\ref{cxb25}), and Gronwall's inequality yields that
\be  \la{cxb67}\ba
\sup_{0\le t\le T}\si\|\n^{1/2}u_{tt}\|_{L^2}^2+\int_{0}^T\si^2\|\nabla u_{tt}\|_{L^2}^2dt\le C.
\ea\ee
Furthermore, it follows from \eqref{cxb480} and \eqref{cxb25} that
\be  \la{cxb68}\ba
\sup_{0\le t\le T}\si\|\nabla u_t\|_{H^1}^2\le C.
\ea\ee

Finally, we deduce from \eqref{cxb52}, \eqref{cxb53}, \eqref{cxb25}, \eqref{cxb44}, \eqref{cxb45}, \eqref{cxb67} and \eqref{cxb68}
 that
\bnn  \la{cxb69}\ba
\si\|\na^2 u\|_{W^{1,q}}
& \le C(\si +\si\|\na  u_t\|_{L^2}+\si\| \na(\n
\dot u )\|_{L^{q}}+\si\|\na^2  P\|_{L^{q}})\\
& \le C(\sigma+\sigma^{\frac{1}{2}} +  \si\|\na u\|_{H^2}+\sigma^{\frac{1}{2}}(\sigma\|\na u_t\|_{H^1}^2)^{\frac{3(q-2)}{4q}})\\
&\le C\sigma^{\frac{1}{2}}+C\sigma^{\frac{1}{2}}(\sigma^{-1})^{\frac{3(q-2)}{4q}}\\
&\le C ,
\ea\enn which
together with (\ref{cxb67}) and (\ref{cxb68}) leads to (\ref{cxb58}) and completes the proof of Lemma \ref{xle5}.
\end{proof}

\section{\la{se6}Proofs of  Theorems  \ref{th1}-\ref{th3}}

With all the a priori estimates in Sections \ref{se3} and   \ref{se5} at hand, we are going to  prove the main results, Theorems \ref{th1}--\ref{th3}.

  \subsection{Proof of Theorem \ref{th1} }  By Lemma \ref{loc1}, there exists a
$T_*>0$ such that the  system (\ref{a1})-(\ref{ch1}) has a unique classical solution $(\rho,u)$ on $\Omega\times
(0,T_*]$. One may use the a priori estimates, Proposition \ref{pr1} and Lemmas \ref{xle3}-\ref{xle5} to extend the classical
solution $(\rho,u)$ globally in time.

First, by the definition of $A_1(T)$, $A_2(T)$ (see \eqref{As1}, \eqref{As2}), the assumption of the initial data \eqref{dt2} and \eqref{bcbh2}, one immediately checks that
$$ A_1(0)+A_2(0)=0, \,\, 0\leq\rho_0\leq \hat{\rho},\,\, A_3(0)\leq C_0^{\delta_0}.$$
Therefore, there exists a
$T_1\in(0,T_*]$ such that
\be\la{dlbh1}\ba
0\leq\rho_0\leq2\hat{\rho},\,\,A_1(T)+A_2(T)\leq 2C_0^{\frac{1}{3}}, \,\, A_3(\sigma(T))\leq 2C_0^{\delta_0}
\ea\ee
hold for $T=T_1.$

Next, we set
\bnn \la{dlbh2}
T^*\triangleq\sup\{T\,|\,{\rm (\ref{dlbh1}) \ holds}\}.
\enn
Then $T^*\geq T_1>0$. Hence, for any $0<\tau<T\leq T^*$
with $T$ finite, it follows from Lemmas \ref{xle3}-\ref{xle5}
that
 \bnn \la{dlbh3}
   \rho  \in C([0,T]; W^{2,q}), \, \na u_t \in C([\tau ,T]; L^q),\quad
 \na u,\na^2u \in C\left([\tau ,T];
 C (\bar{\Omega})\right), \enn
 where one has taken advantage of  the standard
embedding
$$L^\infty(\tau ,T;H^1)\cap H^1(\tau ,T;H^{-1})\hookrightarrow
C\left([\tau ,T];L^q\right),\quad\mbox{ for any } q\in [1,6).  $$
  This in particular  yields
\be\la{dlbh4} \n^{1/2}u_t, \quad\n^{1/2}\dot u \in C([\tau,T];L^2).\ee
Finally, we claim that \bnn \la{dlbh5}T^*=\infty.\enn Otherwise,
$T^*<\infty$. Then by Proposition \ref{pr1}, it holds that
\bnn\la{dlbh6}\ba
0\leq\rho\leq\frac{7}{4}\hat{\rho} ,\,\,\,A_1(T^*)+A_2(T^*)\leq C_0^{\frac{1}{3}},\,\,\, A_3(\sigma(T^*))\leq C_0^{\delta_0} .
\ea\enn
It follows from Lemmas \ref{xle4}, \ref{xle5} and
(\ref{dlbh4}) that $(\n(x,T^*),u(x,T^*))$ satisfy
the initial data condition (\ref{dt1})-(\ref{dt3}) except $ u(\cdot,T^*)\in H^s,$
where  $g(x)\triangleq\n^{1/2}\dot u(x, T^*),\,\,x\in \Omega.$ Thus, Lemma
\ref{loc1} implies that there exists some $T^{**}>T^*$ such that
(\ref{dlbh1}) holds for $T=T^{**}$, which contradicts the definition of $ T^*.$ 

By Lemmas \ref{loc1} and \ref{xle3}-\ref{xle5}, it indicates that $(\rho,u)$ is in fact the unique classical solution defined on $\Omega\times(0,T]$ for any  $0<T<T^*=\infty.$  Moreover, Proposition  \ref{beha1} gives \eqref{qa1w} and we finish the proof of Theorem \ref{th1}. \thatsall

  \subsection{Proof of Theorem \ref{th2} }
For $T>0$, we introduce the Lagrangian coordinates
  \be \la{c61}  \begin{cases}\frac{\partial}{\partial \tau}X(\tau; t,x) =u(X(\tau; t,x),\tau),\,\,\,\, 0\leq \tau\leq T\\
 X(t;t,x)=x, \,\,\,\, 0\leq t\leq T,\,x\in\bar{\Omega}.\end{cases}\ee
 By virtue of \eqref{dt6}, the transformation \eqref{c61} is well-defined. Therefore, by $\eqref{a1}_1$, we get
 \be\la{c62}\ba
\rho(x,t)=\rho_0(X(0; t, x)) \exp \{-\int_0^t\div u(X(\tau;t, x),\tau)d\tau\}.
\ea \ee
  Since $\n_0(x_0)=0$ for some point $x_0\in \Omega,$  for any $t>0,$ there is a point $x_0(t)\in \bar{\Omega}$ such that $X(0; t, x_0(t))=x_0$. Hence, by \eqref{c62}, $\rho(x_0(t),t)\equiv 0$ for any $t\geq 0.$ As a result of Gagliardo-Nirenberg's inequality \eqref{g2}, we get that for  $ r_1\in  (3,\infty),$
\bnn\la{c63}\ba\bar{\rho_0}\equiv\bar\n\leq\|\rho-\bar{\rho}\|_{C\left(\ol{\O }\right)} \le C
\|\rho-\bar{\rho}\|_{L^2}^{\theta_1}\|\na \rho\|_{L^{r_1}}^{1-\theta_1}
\ea\enn
where $\theta_1=2(r_1-3)/(5r_1-6)$. Combining this with \eqref{qa1w} gives \eqref{qa2w} and completes the proof of Theorem \ref{th2}. \thatsall

  \subsection{Proof of Theorem \ref{th3} }

Let $(\rho_0, u_0)$ be the initial data as in Theorem \ref{th3}.  We construct an approximation initial value $(\rho^\delta_0, u^\delta_0)$ satisfying \eqref{dt7} and for any $p\geq 1$,
\begin{equation*}
\lim_{\delta\rightarrow 0}(\|\rho^\delta_0-\rho_0\|_{L^{p}}+\|u^\delta_0-u_0\|_{H^1})=0,
\end{equation*}
\begin{equation*}
\rho^\delta_0\rightarrow\rho_0\  \mathrm{in}\ W^*\ \mathrm{topology}\ \mathrm{of}\ L^\infty\ \mathrm{as}\ \delta\rightarrow 0.
\end{equation*}
Following the proofs in Section 3, one can check that
\be\la{Add1}\ba
  &\sup_{   0\le t\le T  }\left(\|\rho^\delta\|_{L^\infty}+\|\nabla u^\delta\|_{L^2}^2+\sigma^3\|\rho^\delta u_t^\delta\|_{L^2}\right)\\& + \int_0^{T} \sigma\|\rho^\delta u_t^\delta\|_{L^2}^2dt+ \int_0^{T}\sigma^3\|\nabla u_t^\delta\|_{L^2}^2dt\leq C,
  \ea\ee
and then, by Lemma \ref{le3},
\be\la{Add2}
  \sup_{   0\le t\le T  }(\sigma^3\|\nabla u^\delta\|_{L^6}) + \int_0^{T} \sigma^3\|u_t^\delta\|_{L^6}^2dt\leq C,
  \ee
\be\la{Add3}\ba&
  \sup_{   0\le t\le T  }\left(\sigma^3\|\curl u^\delta\|_{H^1}+\sigma^3\|F^\delta\|_{H^1}\right) \\&+ \int_0^{T} \sigma^3(\|\curl u^\delta\|_{W^{1,6}}^2+\|F^\delta\|_{W^{1,6}}^2)dt\leq C,
 \ea \ee
where $F^\delta\triangleq(\lambda+2\mu)\div u^\delta-P(\rho^\delta)+\overline{P(\rho^\delta)}$.
By virtue of Aubin-Lions Lemma,  there exists a subsequence that
\begin{equation*}
\begin{split}
&u^\delta\rightharpoonup u \mbox{ weakly  * } \mathrm{in}\ L^\infty(0, T;H^1), \\
&u^\delta\rightarrow u\  \mathrm{in}\ C([\tau, T]; L^6),\\
&F^\delta\rightharpoonup F \mbox{ weakly  * }  \ \mathrm{ in }\ L^2(0, T;W^{1,6}),\\
&\omega^\delta\rightarrow \omega, \ F^\delta\rightarrow F\ \mathrm{in}\ L^{2}(\tau, T;L^6),
\end{split}
\end{equation*}
for any  $\tau\in (0, T)$. By  standard arguments (see \cite{Ho3} or \cite{L1}), one can deduce the strong convergence of $\rho^\delta$, that is,
\begin{equation*}
\rho^\delta\rightarrow \rho\ \mathrm{ in }\ C([0,  T];L^q(\Omega)),
\end{equation*}
for any $q\in [1,\infty)$. Moreover, \eqref{111} is established directly.
Therefore, we conclude that $(\rho, u)$ is a   weak solution as in Theorem \ref{th3} which  finishes our proof of Theorem \ref{th3}. \thatsall
\section{\la{se7}Proof of  Theorem  \ref{th4}}
This section is devoted to the proof of Theorem  \ref{th4}. Compared with the previous theorems, since the domain is not necessarily simply connected, the inequality \eqref{paz1} with $k=0$ and $p=2$ is no longer valid. Therefore, the difficulty of this proof is that we need an alternative inequality. Thanks to \cite[Proposition 3.7]{aacg}, the following lemma gives   equivalent  norms of $H^1$.
\begin{lemma}  \la{Nslip1}
Let $\Omega$ be a bounded Lipschitz domain in $\r^3$. Then for $v\in H^1$ with $v\cdot n=0$ on $\partial\Omega$, we have the following equivalence of norms:
\be\la{equiv1}
\|v\|_{H^1}\simeq \begin{cases}
\|D(v)\|_{L^2}, \,\, &\Omega \,\,\text{is not axially symmetric}, \\
\|D(v)\|_{L^2}+\int_{\partial\Omega} u\cdot B\cdot uds, \,\, &\Omega\,\, \text{is axially symmetric},
\end{cases} \ee
where $\simeq$ denotes the equivalence of two norms and $D(v) = (\nabla v+(\nabla v)^{\rm tr})/2$ and $B\in W^{2,6}(\Omega)$ is a positive semi-definite $3\times 3$ symmetric matrix satisfying $B> 0$ on some $\Sigma\subset\partial\Omega$ with $|\Sigma|>0$.
\end{lemma}
\begin{remark}
Similar to what have done in \cite{aacg}, when $\Omega$ is axially symmetric with respect to a constant vector $b\in \r^3$,    we have to add the term $\int_{\partial\Omega} u\cdot B\cdot uds$ in order to exclude such a special case that $v=C\, b\times x$.
\end{remark}

The following lemma enables us to replace the inequality \eqref{paz1} ($k=0$ and $p=2$) with \eqref{Nslip1} so that the previous proofs are still available.
\begin{lemma}  \la{Nslip2}
Let $\Omega$ be a smooth bounded  domain in $\r^3$. Then for $v\in H^2$ with $v\cdot n=0$ on $\partial\Omega$, it holds that
\be\ba\la{equiv2}  2\int D(v)\cdot D(v)dx=2\int(\div v)^{2}dx + \int|\curl v|^{2}dx-2\int_{\partial\Omega} v\cdot D(n)\cdot vds.
\ea\ee
\end{lemma}

\begin{proof} Observe that
$$\Delta v= \nabla\div v -\nabla\times \curl v=2\div(D(v))-\nabla\div v,$$
 which together with \eqref{ch7} gives  \eqref{equiv2}.
\end{proof}

{\bf Proof of Theorem \ref{th4}}. First, it is sufficient to find out where the inequality \eqref{paz1} with $k=0$ and $p=2$ is used, that is, \eqref{eeq1},  \eqref{I4}, \eqref{ax40}, \eqref{xbh9}, \eqref{xbh13}, the proof of
Proposition \ref{beha1} and \eqref{cxb42}.

Then,  we will take advantage of the above two lemmas to deduce the similar  results step by step.

Indeed, setting $B=A+2D(n)$, by Lemma \ref{Nslip1} and the extra assumptions of Theorem \ref{th4}, for any $v\in H^1$ with $v\cdot n=0$ on $\partial\Omega$, we have
\be\la{sub1}
\|\nabla v\|_{L^{2}}^{2}\leq C\left(2\mu\|D(v)\|_{L^{2}}^{2}+\lambda\|\div v\|_{L^{2}}^{2}+\int_{\partial\Omega} v\cdot B\cdot vds\right).
\ee
For \eqref{eeq1},  since $A=B-2D(n)$, by Lemma \ref{Nslip2}, $\phi$ can be rewritten as
\be\ba\la{Neeq1}  \phi=2\mu\|D(u)\|_{L^{2}}^{2}+\lambda\|\div u\|_{L^{2}}^{2}+\int_{\partial\Omega} u\cdot B\cdot uds,
\ea\ee
which, together with \eqref{sub1} gives \eqref{a16}.

Next, \eqref{I4}, \eqref{xbh9}, \eqref{xbh13} and \eqref{cxb42} can be similarly dealt with  to get the results of their next step respectively, and the proof of Proposition \ref{beha1} remains valid if we use \eqref{sub1} instead of \eqref{tdu1}.

Finally, it remains to handle \eqref{ax40}. Setting $v=\dot{u}+(u\cdot\nabla n)\times u^{\perp}$ in \eqref{equiv2},  by \eqref{paz2} and Young's inequality, we deduce from \eqref{ax401}  that
\be\la{Nax40}\ba
&\left(\frac{\sigma^{m}}{2}\|\rho^{\frac{1}{2}}\dot{u}\|_{L^2}^2\right)_t+2\mu\sigma^{m}\|D(v)\|_{L^2}^2+\lambda\sigma^{m}\|\div v\|_{L^2}^2
+\mu\sigma^{m}\int_{\partial\Omega}v\cdot B\cdot vds\\
& \leq Cm\sigma^{m-1}\sigma'(\|\rho^{\frac{1}{2}}\dot{u}\|_{L^2}^2+\|\nabla u\|_{L^2}^2+\|\nabla u\|_{L^2}^4)-\left(\int_{\partial\Omega}\sigma^{m}(u\cdot\nabla n\cdot u)Fds\right)_t \\
&\quad+2\delta\sigma^{m}\|\nabla\dot{u}\|_{L^2}^2+C\sigma^{m}\|\rho^{\frac{1}{2}}\dot{u}\|_{L^2}^2(\|\nabla u\|_{L^2}^4+1)\\
& \quad+C(\delta)\sigma^{m}(\|\nabla u\|_{L^2}^2+\|\nabla u\|_{L^2}^6+\|\nabla u\|_{L^4}^4).
\ea\ee
On the other hand, by \eqref{sub1}, we find that
\bnn\ba
\|\nabla \dot{u}\|_{L^{2}}^{2}\leq & C\left(2\mu\|D(v)\|_{L^{2}}^{2}+\lambda\|\div v\|_{L^{2}}^{2}+\int_{\partial\Omega} v\cdot B\cdot vds\right)\\&+C(\|\nabla u\|_{L^2}^2 +\|\nabla u\|_{L^4}^4).\ea
\enn

Combining this with \eqref{Nax40}, we get \eqref{ax401} by letting $\delta$ is suitably small  and thus finish the proof of Theorem \ref{th4}. \thatsall

\section*{Acknowledgements}  The research  is
partially supported by the National Center for Mathematics and Interdisciplinary Sciences, CAS,
National Natural
Science Foundation of China Grant Nos.  11688101,  11525106, 12071200, 11971401, 11871410, and 11871408, and Double-Thousand Plan of Jiangxi Province (No. jxsq2019101008).

\end{document}